\documentclass[10pt]{amsart}

\usepackage[all]{xy}
\usepackage{amsmath,amsthm,amsfonts}
\newtheorem{thm}{Theorem}[section]
\newtheorem{lemma}[thm]{Lemma}
\newtheorem{prop}[thm]{Proposition}
\newtheorem{cor}[thm]{Corollary}

\theoremstyle{definition}
\newtheorem{ex}[thm]{Example}

\theoremstyle{definition}
\newtheorem*{defn}{Definition}

\newcommand{\bears}{\begin{eqnarray*}}
\newcommand{\eears}{\end{eqnarray*}}
\newcommand{\spec}[1]{\ensuremath{\textnormal{Spec}\, #1}}

\newcommand{\Ch}[1]{\ensuremath{\textnormal{Ch}(#1)}}
\newcommand{\CH}[1]{\ensuremath{\textnormal{CH}(#1)}}
\newcommand{\Chk}[2]{\ensuremath{\textnormal{Ch}_{#1}(#2)}}
\newcommand{\CHk}[2]{\ensuremath{\textnormal{CH}_{#1}(#2)}}

\newcommand{\ChG}{\ensuremath{\textnormal{Ch}^G}}

\newcommand{\one}{\ensuremath{\mathfrak{1}}}
\newcommand{\PP}{\ensuremath{\mathbb{P}}}
\newcommand{\ZZ}{\ensuremath{\mathbb{Z}}}
\newcommand{\Zps}{\ensuremath{(\ZZ /p\ZZ )^{\times}}}

\newcommand{\NN}{\ensuremath{\mathbb{N}}}
\newcommand{\Rp}{\ensuremath{R_0\oplus R_1\oplus\cdots\oplus R_{p-1}}}
\newcommand{\Gm}{\ensuremath{\mathbb{G}_m}}
\newcommand{\Aff}{\ensuremath{\mathbb{A}^1}}

\newcommand{\bu}[1]{\ensuremath{\overline{U}_{#1}}}
\newcommand{\km}[1]{\ensuremath{F(#1)^{\times}/F(#1)^{\times p}}}
\newcommand{\wt}[1]{\ensuremath{\widetilde{#1}}}

\newcommand{\smalliff}{\ensuremath{\Leftrightarrow}}

\SelectTips{cm}{}

\begin{document}

\title{A New Definition of the Steenrod Operations in Algebraic Geometry}
\author{Alex Boisvert}
\date{}
\address{Department of Mathematics, University of California,
        Los Angeles, CA 90095-1555} \email {boisvert@ucla.edu}
\begin{abstract}
The Steenrod operations (mod $p$) in Chow theory are defined for any prime $p$ for a quasi-projective scheme, without appealing to the results of any domain but Milnor's $K$-theory.  The new definition also gives a direct formula that depends only on the scheme itself.  Additionally, basic properties of the operations are proved from the new definition.  The idea is based on a construction of M. Rost.
\end{abstract}

\maketitle

\section{Introduction}\label{intro}

In topology, the Steenrod operations were defined in \cite{steenrod} for singular cohomology of topological spaces.  For a long time, it was not clear that there was an analog in algebraic geometry of singular cohomology, much less the Steenrod operations.  However, these questions were resolved when V. Voevodsky defined these operations in motivic cohomology for use in his proof of the Milnor conjecture.  The Chow groups, a special case of motivic cohomology, were left intact by the Steenrod operations; however, at this time there was no elementary definition of the operations on the Chow groups.  This was taken care of later, when P. Brosnan defined these operations specifically for Chow groups, avoiding the machinery of motivic cohomology.

However, Brosnan's construction is not without its drawbacks.  While the Steenrod operations he defines are defined on the Chow groups, he uses the field of equivariant Chow groups in his definition.  Furthermore, his construction involves embedding a scheme into a smooth space and then showing that the construction is independent of the embedding.  Of course, neither is a tremendous problem, but ideally, one would like to avoid such concerns in any definition and give a direct formula.

More recently, A. Merkurjev has defined the Steenrod operations (mod 2) in a more elementary way, relying only on Chow theory.  The results are more difficult to come by mod $p$ (where $p$ is prime) for a variety of reasons; for example, there may not be $p$-th roots of unity in the ground field.  This paper will define these operations for any prime $p$ over any field of characteristic not equal to $p$.

\subsection{Outline}We fix a prime $p$ and a ground field $F$ with char $F\neq p$.  In section \ref{mupactions}, the action of $G=\mu_p$ on schemes is discussed in an elementary way, in the hope that the reader who has not studied algebraic groups in full generality will understand.  Several examples are given, including the action of $G$ on $R^p(X)$ (Example \ref{RpX}) for a scheme $X$.  This action will allow us to define the Steenrod operations over a field without $p$-th roots of unity.  The action of $G$ on a deformation variety is discussed in \ref{defvar}, which later will allow us to reduce to the case of cones (Corollary \ref{reductionvectorbundles}).

In section \ref{torsors}, the concept of a $G$-torsor is introduced, again in an elementary way.  This allows us to introduce the homomorphism $\alpha_T:\CH{T/G}/p\to \CH{T/G}/p$ for any $G$-torsor $T\to T/G$; this map essentially arises from the fact that $G$-torsors over a field $L$ are classified by $L^{\times}/L^{\times p}$.  In section \ref{rostop}, the Rost operation is introduced; it is the building block of the Steenrod operations.  The construction is essentially as follows: if $W$ is a variety over $F$ with a $G$-action, then $W\setminus W^G\to (W\setminus W^G)/G$ is a $G$-torsor; we set $\bu{W}=(W\setminus W^G)/G$.  The Rost operation $\rho_i^W$ is then defined as the composition $\partial^{\bu{W}}_{X}\circ e(L)^{i-1}\circ\alpha_{\bu{W}}:\CHk{*}{\bu{W}}/p\to\CHk{*-i}{W^G}/p$, where $L$ is the canonical line bundle over $\bu{W}$.  Certain properties of the Rost operation are discussed, including functorialities.

In section \ref{steenroddef}, a map $P_X:\CH{X}\to\CH{(R^p(X)\setminus X)/G}$ is defined that essentially mimics Brosnan's map $P^n_G$ from \cite{brosnan}.  At this point we can define the Steenrod operations as follows:
\[S_k^X=(-1)^{n-k}\rho_{(p-1)(n+k)}\circ P_X:\CHk{n}{X}/p\to\CHk{n-(p-1)k}{X}\]
where $\rho=\rho^{R^p(X)}$ (this is possible as $R^p(X)^G=X$, a fact proved in Example \ref{RpX}).

Section \ref{compsteenrod} gives an explicit calculation of the Steenrod operations, which proves that these operations are the same as Brosnan's.  In section \ref{steenrodprop} some properties are proved; finally, section \ref{smoothsteenrod} defines the Steenrod operations of cohomological type and proves further properties.

\section{$\mu_p$-actions on algebraic schemes}\label{mupactions}

In this section we define a $\mu_p$-action on a scheme, give some examples and prove some properties.

\subsection{Definitions and basic properties}

Fix a prime $p$.  Let $X=\spec{R}$ be an affine scheme over a field $F$.
Suppose that $R$ is a \(\ZZ /p\ZZ\)-graded ring \(R=R_0\oplus
R_1\oplus\cdots\oplus R_{p-1}\).  Then the group
\(G=\mu_p=\spec{F[t]/(t^p-1)}\) acts on $X$:
the co-action $F$-algebra homomorphism
\[\theta :\Rp\to F[t]/(t^p-1)\otimes (\Rp )\]
is defined by the formula
\[\theta (r_0+r_1+\cdots +r_{p-1})=1\otimes r_0+t\otimes r_1+\cdots
+t^{p-1}\otimes r_{p-1}\]

Conversely, a $G$-action on $X$ is induced by a \(\ZZ /p\ZZ\)-graded ring
on $R$ as above.  We call such a scheme an \textit{affine $G$-scheme}.
Note that the action is trivial if and only if \(R_1\oplus
R_2\oplus\cdots\oplus R_{p-1}=0\).

In general, a {\itshape $G$-scheme} is a scheme $X$ together with a $G$-action
on $X$.  If $X$ is a quasi-projective scheme, every pair of points on $X$
belongs to some open affine subscheme.  It follows that there is a
$G$-invariant affine covering.

From now on we only consider quasi-projective $G$-schemes over a field.
In the constructions and proofs we may restrict to the class of affine
$G$-schemes.

\begin{defn} Let $X$ be an affine $G$-scheme; consider the ideal
\[I=\sum_{i+j=p}R_iR_j\subset R_0\]
We call the scheme $\spec{R_0/I}$ the \textit{fixed-point subscheme} of $X$,
denoted $X^G$.  Notice that 
\[R_0/I=R/\!\! <\!\! R_1\oplus
R_2\oplus\cdots\oplus R_{p-1}\!\! >\, =R/(I\oplus 
R_1\oplus\cdots\oplus R_{p-1})\]
The natural morphism $\phi :X^G\to X$
satisfies the following universal property: if $Y$ is an affine scheme
with trivial $G$-action, then every $G$-equivariant morphism $Y\to X$
factors uniquely through $\phi$.
\end{defn}

\begin{defn} Again let $X$ be an affine $G$-scheme.  Denote the scheme
$\spec{R_0}$ by $X/G$, called the \textit{quotient} (of $X$ by $G$).  The
natural morphism $\pi :X\to X/G$ satisfies the following universal
property: if $Y$ is an affine scheme with trivial  $G$-action, then every
$G$-equivariant morphism $X\to Y$ factors uniquely through $\pi$.
\end{defn}

\subsection{Examples}

\begin{ex}\label{RpX}
\textit{The scheme $\mathit{R^p(X)}$.}  Consider the 
algebra \(F_p=F[x]/(x^p-1)\).  For any scheme $X$ over $F$, we write 
$R^p(X)$ for the scheme \(R_{F_p/F}(X_{F_p})\), representing the functor 
\(S\mapsto X(S\otimes_F F_p)\).  Note that if char$(F)\neq p$ and $F$ 
contains all $p$-th roots of unity, then $R^p(X)$ is a product of $p$ 
copies of $X$.  

We have a canonical $G$-action on $R^p(X)$ defined as follows: to define
\(G\times R^p(X)\to R^p(X)\), we can define a map \(G(S)\times
R^p(X)(S)\to R^p(X)(S)\), i.e.\ \(\mu_p(S)\times X(S\otimes F_p)\to
X(S\otimes F_p)\).  For $\xi\in\mu_p(S)$ we can define
\(\gamma_{\xi}:S\otimes F_p\to S\otimes F_p\) via $a\otimes x^m\mapsto
a\xi^m\otimes x^m$.  Hence our map is of the following form:
\bears
\mu_p(S)\times X(S\otimes F_p) &\to& X(S\otimes F_p) \\
(\xi ,f) &\mapsto & f\circ (\gamma_{\xi})^{\#}
\eears
Suppose  char$(F)\neq p$ and $F$ contains a $p$-th root of unity.  Then 
this action is just the standard action of permuting coordinates on $X^p$, 
provided we have selected a distinguished $p$-th root of unity.  

By the generalized Yoneda Lemma (\cite{eisenbud} Proposition VI-2), 
the embedding \(X(S)\hookrightarrow X(S\otimes F_p)\) gives rise to a 
morphism $X\to R^p(X)$.  If char$(F)\neq p$ and $F$ contains a $p$-th 
root of unity, this morphism corresponds with the diagonal morphism 
$X\to X^p$.  Hence, if char$(F)\neq p$, by descent theory, we
see that $X$ is identified with $(R^p(X))^G$.
\end{ex}

\begin{ex}\label{cone}
Let $C=\spec{S^{\bullet}}$ be a cone over $X=\spec{S^0}$.  For each 
$i=0,1,\ldots ,p-1$, let $R_i$ be the coproduct of all $S^j$ with $j\equiv 
i$ (mod $p$).  The $\ZZ /p\ZZ$-grading on $S^{\bullet}=\Rp$ gives rise to 
a $G$-action on $C$ which we call the \textit{canonical $G$-action.}

By definition, $C^G=\spec{R_0/I}$.  Since $S^i=(S^1)^i$, it follows that 
$C^G=\spec{S^0}$, the image of the zero section of $C$.  We also have that 
$C/G=\spec{R_0}$.  If $C$ is a line bundle over $X$, i.e.\ $S^1$ is an 
invertible sheaf and $S^i=(S^1)^{\otimes i}$, then $C/G=C^{\otimes p}$.
\end{ex}

\subsection{The deformation variety}\label{defvar}

Let $X=\spec{R}$ be an affine scheme with a $G$-action; let 
$Y=X^G=\spec{(R/J)}$ where $J=I\oplus R_1\oplus\cdots\oplus R_{p-1}$ as 
before.  Let $N$ be the normal cone of the closed embedding $Y\to X$ and 
let
\[D=\spec{\wt{R}}=\spec{(\cdots\oplus J^2t^{-2}\oplus Jt^{-1}\oplus 
R\oplus Rt\oplus\cdots )}\]
be the deformation variety.  We have a canonical $\ZZ /p\ZZ$-grading on 
$\wt{R}$ induced by the grading on $R$, giving $t$ the zero grading.  
Explicitly, for each $i=0,1,\ldots p-1$, we have
\[\wt{R}_i=\cdots\oplus (J^2)_it^{-2}\oplus J_it^{-1}\oplus
R_i\oplus R_it\oplus\cdots \]

We consider the ideal \(\wt{I}=\sum_{i+j=p}\wt{R}_i\wt{R}_j\subset 
\wt{R}_0\).  For $n\geq 0$, $\wt{I}_n=I$.  Hence we turn our attention 
to $\wt{I}_{-k}$ for $k>0$.

\begin{prop}\label{defcalc} $\wt{I}_{-k}=(J^k)_0$ for each $k>0$.  
\end{prop}
\begin{proof}
One implication is clear: $\wt{I}_{-k}\subset (J^k)_0$.  It is also 
clear that 
\[\sum_{i+j=p}J_i(J^{k-1})_j\subset \wt{I}_{-k}\]
so it is enough to show that $(J^k)_0\subset \sum_{i+j=p}J_i(J^{k-1})_j$.  
Now 
\[(J^k)_0=(J\cdot J^{k-1})_0=J_0\cdot 
(J^{k-1})_0+\sum_{i+j=p}J_i(J^{k-1})_j\]
so we are reduced to showing that
\(J_0\cdot (J^{k-1})_0\subset \sum_{i+j=p}J_i(J^{k-1})_j\).  
We proceed by induction.  For $k=1$, the left-hand side becomes 
\(J_0\cdot (J^0)_0=IR_0=I\)
and the right-hand side becomes 
\[\sum_{i+j=p}J_i(J^0)_j=\sum_{i+j=p}R_iR_j=I\]
so the base step is complete.

Now assume that $J_0\cdot (J^{k-1})_0\subset \sum_{i+j=p}J_i(J^{k-1})_j$.  
Then 
{\allowdisplaybreaks
\bears
J_0(J^k)_0 &=& J_0(J\cdot J^{k-1})_0 \\
&=& J_0\cdot\left[ J_0\cdot (J^{k-1})_0+\sum_{i+j=p}J_i(J^{k-1})_j\right] \\
&\subset & J_0\cdot\left[ \sum_{i+j=p}J_i(J^{k-1})_j\right] \subset \sum_{i+j=p}J_i(J^{k})_j
\eears
}
completing the proof.
\end{proof}
\begin{cor}\label{deformationaction} If $D$ is the deformation variety as 
above, then $D^G=X^G\times\Aff$.
\end{cor}
\begin{proof}
By definition, $D^G=\spec{(\wt{R}_0/\wt{I})}$.  But we have that
\[\wt{R}_0=\cdots\oplus (J^2)_0t^{-2}\oplus J_0t^{-1}\oplus
R_0\oplus R_0t\oplus\cdots \]
and, by Proposition \ref{defcalc},
\[\wt{I}=\cdots\oplus (J^2)_0t^{-2}\oplus J_0t^{-1}\oplus
I\oplus It\oplus\cdots \]
so
\[\wt{R}_0/\wt{I}=(R_0/I)\oplus (R_0/I)t\oplus (R_0/I)t^2\oplus\cdots\]
giving us that $D^G=\spec(\wt{R}_0/\wt{I})=X^G\times\Aff$.\end{proof}

\section{$G$-torsors and the homomorphism $\mathbf{\alpha}$}\label{torsors}

Now that we have introduced the concept of a group action, we move on to $G$-torsors, the ``nice'' group actions.  We define $G$-torsors, give some examples and prove some basic properties.  Then we can define, for any $T\to T/G$ a $G$-torsor, the homomorphism $\alpha_T$ which is essential to our main construction.

\subsection{Notation and definition}

\noindent \textbf{Notation.}  Let $K_*(F)$ be the graded Milnor ring of a 
field.  For a scheme $X/F$, we set
\[C(X)=\coprod_{x\in X}K_*F(x)\]
C(X) is a group, graded by the dimension of $x\in X$.  We also set
\[C_{q,r}(X)=\coprod_{x\in X_{(q)}}K_{q+r}F(x).\]
The $q$-th homology group of the complex $C_{q,r}(X)$ is denoted by $A_q(X,K_r)$; 
these are called the \textit{$K$-homology groups}.  Finally, we define 
$\textrm{CH}_n(X)=A_n(X,K_{-n})$ and $\CH{X}=\coprod_n\textrm{CH}_n(X)$, the 
\textit{Chow groups of X}.

We will also fix a prime $p$ in this paper and simplify notation as follows:
for $\textrm{CH}_n(X)/p$ and $\CH{X}/p$ we will write $\textrm{Ch}_n(X)$ and 
$\Ch{X}$ respectively.  Furthermore, instead of $A_p(X,K_r/p)$, we will write 
$A_{q,r}(X)$.

\begin{lemma}\label{tensorsurjective} Let $R$ be a ring and let $f:M\to N$ 
be a morphism of projective $R$-modules.  Suppose $f^{\otimes 
k}:M^{\otimes k}\to N^{\otimes k}$ is surjective for some $k\geq 1$.  Then 
$f$ itself is surjective.
\end{lemma}
\begin{proof} Let $\mathfrak{m}\subset R$ be a maximal ideal.  When we 
localize at $\mathfrak{m}, M_{\mathfrak{m}}^{\otimes k}\to 
N_{\mathfrak{m}}^{\otimes k}$ is now a surjective map of a tensor power of 
free modules.  Hence $M_{\mathfrak{m}}\to N_{\mathfrak{m}}$ must be 
surjective.  Now, since this is true for all maximal ideals 
$\mathfrak{m}$, $M\to N$ must be surjective.
\end{proof}

\begin{prop}\label{torsordef} Let $X=\spec{R}=\spec{(\Rp )}$ be an affine 
$G$-scheme.  Let 
$I=\sum_{i+j=p}R_iR_j\subset R_0$ as before.  Then the following are 
equivalent:

(1) $X^G=\emptyset$

(2) $I=R_0$

(3) $R_i^p=R_0$ for all $i=1,2,\ldots ,p-1$.

(4) $R_1^p=R_0$

(5) $R_i\!\cdot\! R_{p-i}=R_0$ for each $i=1,2,...p-1$.

(6) The product map $R_k\otimes_{R_0}\! R_l\to R_{k+l}$ 
(where the sum is taken mod $p$) is an isomorphism for all 
$k,l=0,1,2,\ldots ,p-1$.

(7) The homomorphism $\Phi: R\otimes_{R_0}\! R\to F[t]/(t^p-1)\otimes R$ 
given by 
\[a_k\otimes b_l\mapsto t^k\otimes a_k\!\cdot\! b_l\,\,\,\,\,(a_k\in R_k, 
b_l\in R_l)\]
is an isomorphism.
\end{prop}

\begin{proof}
Certain implications are clear: $(1)\smalliff (2), (5)\smalliff (6), 
(3)\Rightarrow
(4)\Rightarrow (5)\Rightarrow (2)$, $(6)\Rightarrow (5)$, and 
$(7)\smalliff (6)$.  It will hence suffice to 
show that $(2)\Rightarrow (3)$ and $(3)+(5)\Rightarrow (6)$.

$\mathbf{(2)\Rightarrow (3):}$  The result is obvious if $p=2$, so let $p$ 
be an odd prime.  Choose an index $d$.  We already have $R_d^p\subset R_0$ 
so it suffices to show that 
$R_0\subset R_d^p$.

Let $m=\frac{p-1}{2}$ and $N=p!m$.  Then
\[R_0=R_0^N=(\sum_{i+j=p}R_iR_j)^N=\sum_{\stackrel{n_i\geq 0\,\forall 
i}{\sum n_i=N}}\!\!c_{\alpha}\!\!\prod_{i+j=p}(R_iR_j)^{n_i}\]

for some constants $c_{\alpha}$.  Hence it suffices to show
\[\prod_{i+j=p}(R_iR_j)^{n_i}\subset R_d^p\]
for any choice of $(n_1,\ldots ,n_m)\in\NN^m$ that sum to $N$.

Choose such a multi-index $(n_1,\ldots ,n_m)$.  From our choice of $N$, 
there exists a $k$ such that $n_k\geq p!$.  Hence we have
\[\prod_{i+j=p}(R_iR_j)^{n_i}=\prod_{\stackrel{i+j=p}{i\neq 
k}}(R_iR_j)^{n_i}(R_kR_{p-k})^{n_k}\subset 
(R_kR_{p-k})^{n_k}\subset (R_kR_{p-k})^{p!}\]

There exist $m_1,m_2\in\{ 1,2,\ldots ,p-1\}$ (specifically, the elements 
in $(\ZZ /p\ZZ)^{\times}$ such that $km_1=d$ and $(p-k)m_2=d$) such that 
$(R_k)^{m_1}\subset R_d$ and $(R_{p-k})^{m_2}\subset R_d$.  Thus:
\[(R_kR_{p-k})^{p!}=[ 
(R_k^{m_1})^{\frac{(p-1)!}{m_1}}(R_{p-k}^{m_2})^{\frac{(p-1)!}{m_2}}]^p\subset 
[ (R_d)^{\frac{(p-1)!}{m_1}}(R_d)^{\frac{(p-1)!}{m_2}}]^p\subset R_d^p\]
the desired result.

$\mathbf{(3)+(5)\Rightarrow (6):}$ Take $\sum x_i\otimes y_i$ in the 
kernel of the product map, i.e.\  $\sum x_iy_i=0$.  Choose $a_j\in R_k$ and 
$b_j\in R_{p-k}$ such that $\sum a_jb_j=1$.  Then $b_jx_i\in R_0$ and so 
\[\sum x_i\otimes y_i=\sum a_jb_jx_i\otimes y_i=\sum a_j\otimes b_jx_iy_i=0\]
hence the product map is injective.

Next, note that by slightly modifying the argument above, we can prove 
that the product map $R_i^{\otimes p}\to R_0$ is injective for all $i$, 
and is thus an isomorphism in view of (3).  In particular, each 
$R_i$ is a locally free $R_0$-module of rank 1, hence projective.

Consider the following diagram:

\[
\xymatrix{
(R_k\otimes_{R_0}R_l)^{\otimes p} \ar[r] \ar[d] 
  & (R_{k+l})^{\otimes p} \ar[d] \\
R_0\otimes_{R_0}R_0 \ar[r] & R_0
}
\]

The vertical maps are isomorphisms by our previous discussion, and the 
bottom map is surjective (actually an isomorphism).  Hence the top map is 
a surjection, so by Lemma \ref{tensorsurjective}, $R_k\otimes_{R_0}R_l\to 
R_{k+l}$ is surjective, the desired result.
\end{proof}

\begin{defn}
Let $X$ be a $G$-scheme.  We call $f:X\to X/G$ a \textit{G-torsor} if there exists 
an open affine covering \(X/G=\cup U_i\) such
that \(f^{-1}(U_i)\) satisfies the properties of the previous Proposition.
\end{defn}

\subsection{Examples}

\begin{ex}\label{trivial} If $X$ is any scheme, the scheme $G\times X$ has a 
natural $G$-action: the map $G\times (G\times X)\to G\times X$ is defined 
by $(g_1,g_2,x)\mapsto (g_1g_2,x)$.  The second projection 
$G\times X\to X$ is a $G$-torsor, known as the \textit{trivial $G$-torsor.}
\end{ex}

\begin{ex}\label{UGtorsor} Let $X$ be a $G$-scheme and $U=X\setminus X^G$.  
Since $U^G$ is empty, the morphism $U\to U/G$ is a $G$-torsor.
\end{ex}

\begin{ex}\label{tp-r} Suppose $X$ is an affine $G$-scheme, $X\to X/G$ a 
$G$-torsor, and $R_0$ is a local ring.  Then $R_1$ is a free $R_0$-module 
of rank 1, i.e.\ $R_1=bR_0$ for some $b\in R^{\times}$.  We have that, 
for each $i=1,2,\ldots ,p-1$, $b^i$ generates $R_i$ (as $b^{\otimes i}$ 
generates $R_1^{\otimes i}$).  Let $a=b^p\in R_0^{\times}$.  Then $R$ is 
a free module on $e_0,e_1,\ldots ,e_{p-1}$, where $e_i=b^i$, and
\[
e_ie_j=
\begin{cases}
e_{i+j}, & i+j\leq p-1 \\
ae_{i+j-p}, & i+j\geq p
\end{cases}
\]
that is, $R\cong R_0[t]/(t^p-a)$.  Note that $a$ is unique up to a $p$-th 
power.
\end{ex}

\begin{ex}\label{fibertorsor} Let $X\to X/G$ be a $G$-torsor, $X$ affine.  Let 
\(R_0\to S_0\) be a ring homomorphism and set $S=R\otimes_{R_0}S_0$.  Then 
$(S_1)^p=S_0$ and so the natural morphism \(\spec{S}\to\spec{S_0}\) is a 
$G$-torsor.  In particular, for every point \(y\in X/G\), $X_y$ is a 
$G$-torsor over $\spec{F(y)}$.
\end{ex}

\begin{ex}\label{ff} Let $X\to X/G$ be a $G$-torsor, $X$ affine.  Then each 
$R_i$ is an invertible $R_0$-module, that is, each $R_i$ is locally free 
of rank 1 over $R_0$.  It follows that, in general, $X\to X/G$ is a flat 
morphism.  Since we have in any case that $X\to X/G$ is surjective, $X\to 
X/G$ is in fact faithfully flat.

This implies the following useful result: suppose $X\to X/G$ and $Y\to 
Y/G$ are $G$-torsors; $f:X\to Y$ a morphism.  Let 
$f_G:X/G\to Y/G$ be the induced map on the quotient spaces.  Let 
$\mathcal{P}$ be one of the following properties: proper, flat, smooth, 
regular embedding or l.c.i morphism.  Then if $f$ has property 
$\mathcal{P}$, by descent theory, so does $f_G$.
\end{ex}

\subsection{Basic properties}

\begin{prop}\label{gtorsorcones} Let $T\to T/G$ and $S\to S/G$ be 
$G$-torsors, with $f:S\to T$ a closed embedding.  Then:

(a) The induced map $g:S/G\to T/G$ is a closed embedding, and

(b) let $N_f$ and $N_g$ be the normal cones of $f$ and $g$, respectively.  
Then $N_g=N_f/G$.
\end{prop}
\begin{proof}
We may assume that all the schemes are affine.  Then $T=\spec{R}$, 
$T/G=\spec{R_0}$, $S=\spec{R/K}$ for some (graded) ideal $K\subset R$ and 
$S/G=\spec{R_0/K_0}$.  Hence (a) follows easily.  Note that we did not use 
the assumption that these were torsors.

Further, we have
\bears
N_f &=& \spec{(R/K\oplus K/K^2\oplus K^2/K^3\oplus\cdots )} \\
N_g &=& \spec{(R_0/K_0\oplus K_0/K_0^2\oplus K_0^2/K_0^3\oplus\cdots )} \\
\textrm{and}\,\, N_f/G &=& \spec{(R_0/K_0\oplus (K/K^2)_0\oplus 
(K^2/K^3)_0\oplus\cdots )}
\eears
so for (b) it suffices to show $(K^n)_0=K_0^n$ for each $n\in\NN$.  We 
have that $K_0^n\subset (K^n)_0$, so it remains to show $(K^n)_0\subset 
K_0^n$.

Take $a_1,a_2,\ldots ,a_n\in K$ such that $a=\prod_{i-1}^na_i\in R_0$; $a$ 
is a typical element of $(K^n)_0$.  Write each $a_i$ as 
$a_i=a_i^0+a_i^1+\cdots +a_i^n$ with $a_i^j\in K_j$.  Then
\[a=\prod_{i=1}^na_i=\sum_{p|j_1+\cdots +j_n}a_1^{j_1}a_2^{j_2}\cdots 
a_n^{j_n}\]
the other terms having cancelled.  Hence it is enough to show that 
whenever $p|j_1+j_2+\cdots +j_n$, $K_{j_1}K_{j_2}\cdots K_{j_n}\subset 
K_0^n$.  It follows from Proposition \ref{torsordef} that 
$R_{p-j_1}R_{p-j_2}\cdots R_{p-j_n}=R_0$, so
\bears
K_{j_1}K_{j_2}\cdots K_{j_n} &=& R_0K_{j_1}K_{j_2}\cdots K_{j_n} \\
&=& (R_{p-j_1}K_{j_1})(R_{p-j_2}K_{j_2})\cdots (R_{p-j_n}K_{j_n}) \\
&\subset & K_0K_0\cdots K_0=K_0^n
\eears
\end{proof}

\begin{prop}\label{gtorsorprops}
Suppose char$(F)\neq p$ and F contains a (primitive) $p$-th root of 
unity. 
We denote $G'=G(F)=\mu_p(F)$, the abstract group of the roots of unity of 
$F$.
Let $\pi:X\to X/G$ be a $G$-torsor.  Then:

(1) $\pi^*$ induces an isomorphism $\pi^*:Z(X/G)\to Z(X)^{G'}$.

(2) $\pi_*\circ\pi^*:Z(X/G)\to Z(X/G)$ is multiplication by $p$.

(3) Choose a distinguished generator $\xi$ of $G'$ and let $\sigma$ be 
the automorphism of $Z(X)$ given by $\xi$.  Then 
$\pi^*\circ\pi_*=1+\sigma +\sigma^2+\cdots +\sigma^{p-1}$.
\end{prop}
\begin{proof}
Let $Y=X/G$.  For every point $y\in Y$, \(X_y\to\spec{F(y)}\) is a
$G$-torsor by Example \ref{fibertorsor}. By Example \ref{tp-r}, 
$X_y=\spec{F(y)[t]/(t^p-r)}$ for some $r\in F(y)^{\times}$.  We consider 
two cases.

\textbf{Case 1: $t^p-r$ is irreducible in $F(y)[t]$.}  Then $F(y)[t]/(t^p-r)$ 
is a field and $X_y$ consists of only one point $x$.  Note that 
$\pi^*([y])=[x]$ and $\pi_*([x])=p[y]$.  Furthermore, our automorphism 
$\sigma$ acts on fibers, implying that $\sigma^i([x])=[x]$ for 
$i=0,1,\ldots ,p-1$.  Hence we get (2) and (3) for this case.

\textbf{Case 2: $t^p-r$ is not irreducible in $F(y)[t]$.}  Here $t^p-r$ must 
split completely as $F$ contains all roots of unity, so $F(y)=F\times 
F\times\cdots\times F$ ($p$ times) and $X_y$ has $p$ points 
$x_1,x_2,\ldots ,x_p$.  We see that $\pi^*([y])=[x_1]+[x_2]+\cdots +[x_p]$ 
and, for each $i$, $\pi_*([x_i])=[y]$.  So $\pi_*\circ\pi^*([y])=p[y]$ and 
$\pi^*\circ\pi_*([x_i])=[x_1]+[x_2]+\cdots +[x_p]=(1+\sigma 
+\sigma^2+\cdots +\sigma^{p-1})([x_i])$, as $\sigma$ acts transitively on 
the fiber.  Hence (2) and (3) are proved in this case.

One final thing to note is that cycles of the form $[x]$ and 
$[x_1]+[x_2]+\cdots +[x_p]$ (as above) generate the group $Z(X)^{G'}$ of 
$G'$-invariant cycles.  Hence $\pi^*$ is surjective.  But since it is 
certainly injective, it is an isomorphism.
\end{proof}

Let $X$ be any scheme with a $G$-action. We have two morphisms $G\times 
X\to X$; the action morphism $a$, which (roughly speaking) sends 
$(g,x)\mapsto g\cdot x$ and the second projection $q$, which sends 
$(g,x)\mapsto x$.  If $\pi :X\to X/G$ is the projection, we note that 
$\pi\circ a=\pi\circ q:G\times X\to X/G$, 
as either way, the corresponding $F$-algebra homomorphism $R_0\to 
F[t]/(t^p-1)\otimes R$ sends $r_0\mapsto1\otimes r_0$.

Consider the exact sequence
\[0\to K\to Z(X)\stackrel{a^*-q^*}{\to}Z(G\times X)\]
where $K=\textrm{Ker}(a^*-q^*)\subset Z(X)$.  We denote this kernel as 
$Z(X)^G$.  Note that if $F$ contains a $p$-th root of unity, then
$Z(X)^G=Z(X)^{G'}$, in the notation of Proposition \ref{gtorsorprops}.

\begin{defn} We call a cycle $\delta\in Z(X)$ \textit{G-invariant} if 
$\delta\in Z(X)^G$.
\end{defn}

\begin{cor}\label{torsorisomorphism} Suppose char$(F)\neq p$; let $\pi 
:X\to X/G$ be a 
$G$-torsor.  
Then $\pi^*$ induces an isomorphism $Z(X/G)\to Z(X)^G$.
\end{cor}
\begin{proof} Since $\pi\circ a=\pi\circ q:G\times X\to X/G$, the 
composition 
\[Z(X/G)\stackrel{\pi^*}{\to}Z(X)\stackrel{a^*-q^*}{\to} 
Z(G\times X)\] 
is 0; hence we have a unique map $Z(X/G)\to Z(X)^G$ induced 
by  $\pi^*$.  Now Proposition \ref{gtorsorprops} (1) along with descent 
theory tell us that this map is an isomorphism.
\end{proof}

\subsection{The homomorphism $\alpha$}

Let $T\to T/G$ be a $G$-torsor.  For each \(x\in T/G\) we define
\[a_x\in\km{x}\] in the following way: consider the fiber diagram
\begin{equation}\label{fiberdiagram}
\xymatrix{
\spec {L_x} \ar[r] \ar[d] & T \ar[d] \\
\spec {F(x)} \ar[r] & T/G \\
}
\end{equation}
By Example \ref{tp-r}, \(L_x\cong F(x)[t]/(t^p-a_x)\) for some
uniquely determined $a_x\in\km{x}$.

\begin{defn} For every pair of points \(x,y\in T/G\), we define a
homomorphism
\[\alpha^x_y :K_nF(x)/p\rightarrow K_{n+1}F(y)/p\]
via 
\[
\{c_1,\ldots ,c_n\}\mapsto
\begin{cases}
\{a_y,c_1,\ldots ,c_n\} & \text{if $x=y$} \\
0 & \text{otherwise}
\end{cases}
\]
We hence get an endomorphism \(\alpha_T :C( T/G)\to C( T/G)\) such that
its
$(x,y)$-component is $\alpha^x_y$.  It has degree 0 with respect to the
grading on $C( T/G)$, so we may talk about \(\alpha_T :C_k( T/G)\to
C_k( T/G)\).
\end{defn}

\begin{prop}\label{alphad}
Let $T\to T/G$ be a $G$-torsor.  Then \(d\circ\alpha_T =-\alpha_T\circ d\).
\end{prop}
\begin{proof} 
Consider the diagram
\[
\xymatrix{
C_k( T/G) \ar[r]^-{d} \ar[d]^{\alpha_T} & C_{k-1}( T/G) \ar[d]^-{\alpha_T}
\\
C_k( T/G) \ar[r]^-{d} & C_{k-1}( T/G)
}
\]
The $(x,y)$-component in both directions is nonzero only if $y$ is a
specialization of $x$.  

Let $Z=\overline{\{ x\} }$ and consider the local ring 
$O_{Z,y}$.  It is a local ring, so we may apply our construction and 
get a corresponding $a_{Z,y}\in O_{Z,y}^{\times}$ (unique up to $p$-th power) 
and $L_{Z,y}$ as before.
Note further that $F(x)=$ qf($O_{Z,y}$) and  $F(y)$ is its residue field.
We have that $a_{Z,y}$ specializes to $a_x\in F(x)^{\times}$ and $a_y\in F(y)^{\times}$,
hence the statement follows from \cite{merkurjev}, Proposition 100.4 (1) and Fact 100.8 (3).
\end{proof}

This proposition allows us to consider
\[
\alpha_T :A_m( T/G,K_k/p)\to A_m( T/G,K_{k+1}/p)
\]
or, more succinctly,
\[
\alpha_T :A_{m,k}( T/G)\to A_{m,k+1}( T/G)
\]
When $T$ is clear from context, we will often
use the simpler notation $\alpha$.  
\begin{prop}
\label{alphafunctoriality}Let $g:T'\to T$ be a $G$-equivariant morphism of
$G$-varieties, such that $T\to T/G$ and $T'\to T'/G$ are $G$-torsors.  Let 
$f$ be the induced map $T'/G\to T/G$.  Then:

(1) If $f$ is proper, then \(\alpha_T\circ f_*=f_*\circ\alpha_{T'}\).

(2) If $f$ is flat, then \(\alpha_{T'}\circ f^*=f^*\circ\alpha_{T}\).
\end{prop}

\begin{proof}
Take $x\in T'/G$, and let $y=f(x)\in T/G$.  
We have that $f^*(a_y)=a_x$, so (1) follows from the projection
formula for Milnor's $K$-groups (\cite{merkurjev}, Fact 100.8 (3)) and (2) follows 
immediately from the definitions.
\end{proof}

\begin{prop}\label{alphanaturality}
Let $f:T\to Y=T/G$ be a $G$-torsor and $\alpha_T:A_{m,k}(Y)\to A_{m,k+1}(Y)$ 
be the associated homomorphism.
\begin{enumerate}
\item If $\,T\to Y$ is a trivial $G$-torsor, then $\alpha_T=0$.
\item Suppose $\pi :X\to Y$ is a proper morphism.  Let $S=T\times_YX$ be 
the fiber product.  Then the following diagram is commutative:
\[
\xymatrix{
A_{m,k}(X) \ar[r]^{\alpha_S} \ar[d]^{\pi_*} & A_{m,k+1}(X) \ar[d]^{\pi_*} \\
A_{m,k}(Y) \ar[r]^{\alpha_T} & A_{m,k+1}(Y)
}
\]
\end{enumerate}
\end{prop}
\begin{proof}
1.  Since $T\to Y$ is a trivial torsor, $T=G\times Y$.  Hence in our 
fiber diagram (\ref{fiberdiagram}), \(\spec{L_x}=G\times\spec{F(x)} = 
\spec{(F(x)[t]/(t^p-1))}\), so $a_x=1$ for all $x\in X$.  It follows from 
the definition of $\alpha_T$ that $\alpha_T=0$.

2. Note that $S\to X$ is a $G$-torsor by Example \ref{fibertorsor}.  It 
suffices here to prove the commutativity of the diagram
\[
\xymatrix{
C_l(X) \ar[r]^{\alpha_S} \ar[d]^{\pi_*} & C_l(X) \ar[d]^{\pi_*} \\
C_l(Y) \ar[r]^{\alpha_T} & C_l(Y)
}
\]
The $(x,y)$-component in both directions is nonzero only if $y=\pi (x)$.  Now the statement follows directly from the projection formula for Milnor's $K$-groups.
\end{proof}

\section{The Rost operation}\label{rostop}

The Rost operation (so named because it is a generalization of the construction given by M. Rost in \cite{rost}) builds on $\alpha_T$ defined above and is the main component of our construction of the Steenrod operations.

\subsection{Preliminaries}

Let $F$ a field of characteristic not $p$, \(G=\mu_p\).  Let $W$ be a
variety over $F$ together with a $G$-action. Let \(X=W^G\) be the
fixed-point subscheme, \(U_W=W\setminus X\).

Let \(\overline{W}=W/G\) be the quotient,
\(\overline{U}_W\subset\overline{W}\)
the image of $U$.  We have seen that the projection $\pi :U_W\to\bu{W}$ is a $G$-torsor.  We may consider $X$ as a subscheme of $\overline{W}$, with open complement $\bu{W}$.

\begin{lemma}\label{barf} Suppose $g:S\to W$ is a $G$-equivariant closed embedding
of $G$-schemes.  Let $f:S^G\to W^G$ be the induced morphism.
Then there exists an induced $G$-equivariant closed embedding $\bar{f}:\bu{S}\to\bu{W}$.
Moreover, if $g$ is a regular closed embedding, then so is $\bar{f}$.
\end{lemma}
\begin{proof}
Since $S\cap W^G=S^G$, we have an induced map 
$S\setminus S^G\to W\setminus W^G$.  This in turn induces $\bar{f}:\bu{S}\to\bu{W}$,
which is certainly still a closed embedding.  Moreover,
if $g$ is regular, then so is $S\setminus S^G\to W\setminus W^G$, and so 
$\bar{f}$ will be regular as a result of example \ref{ff}.
\end{proof}

\begin{lemma}\label{opencone} Suppose $f:Z\to X$ is a closed embedding, and 
$U\subset X$ an open subscheme.  We consider $Z$ as a closed subscheme of $X$.
Let $g:Z\cap U\to U$ be the closed embedding and $\pi :N_f\to Z$ be the normal
cone of $Z$ in $X$.  Let $N_g$ be the normal cone of the closed embedding $g$.
Then $N_g=\pi^{-1}(Z\cap U)$.
\end{lemma}
\begin{proof} We may assume that all schemes are affine and that $U\subset X$ is a 
principal open subscheme.  Hence $X=\spec{R}, Z=\spec{R/K}, U=\spec{R_a}$.  We find
$Z\cap U=\spec{R_a/K_a},$
\[N_f=\spec{(R/K\oplus K/K^2\oplus\cdots )}\]
and
\[N_g=\spec{(R_a/K_a\oplus K_a/K_a^2\oplus\cdots )}\]
and so the result is clear.
\end{proof}

\begin{lemma}\label{U_Wcone} Suppose $S\to W$ is a $G$-equivariant closed embedding of
$G$-schemes with normal cone $E$.  Then $\bu{E}$ is the normal cone of the closed embedding
$\bu{S}\to\bu{W}$.\end{lemma}
\begin{proof} We may assume all schemes are affine; i.e., $W=\spec{R}$, $S=\spec{R/K}$, $W^G=\spec{R/J}$, $W^G=\spec{R/(J+K)}$, and 
\[E=\spec{(R/K\oplus K/K^2\oplus K^2/K^3\oplus\cdots )}\]
Now $E^G=\spec{(R/(J+K)\oplus T)}$ for some graded module $T$ that we don't care about; for
the purposes of this proof it is enough to notice that $\pi(E^G)=W^G$.  Hence we can apply
Lemma \ref{opencone} to conclude that $U_E$ is the normal cone of $U_S$ in $U_W$.
Finally, by Proposition \ref{gtorsorcones}, $\bu{E}$ is the normal cone of $\bu{S}$
in $\bu{W}$.
\end{proof}

\subsection{Definition and functoriality}

\begin{defn} Let $L$ be the line bundle on \bu{W} induced from \(\pi
:U_W\to\bu{W}\) via the inclusion \(\mu_p\hookrightarrow\Gm\); that
is, $L=(U_W\times\Aff )/G$ where $G$ acts on $U_W$ as usual and on \Aff\ 
via $G\hookrightarrow\Gm$.  We call $L$ the \textit{canonical line bundle} on $\bu{W}$.

For any \(i\in\NN\), we define the $\mathbf{i}$\textbf{-th Rost operation}
\(\rho^W_i:A_{k,j}(\bu{W})\to
A_{k-i,j+i}(X)\) as the composition
\[
\xymatrix{
A_{k,j}(\bu{W}) \ar[r]^-{\alpha} &
A_{k,j+1}(\bu{W}) \ar[r]^-{e} &
A_{k-i+1,j+i}(\bu{W}) \ar[r]^-{\partial} &
A_{k-i,j+i}(X)
}
\]
where $\alpha =\alpha_{U_W}, e=e(L)^{i-1}$ and $\partial =\partial^{\bu{W}}_X$.
We can define the \textit{total Rost operation $\rho^W$} by taking the
coproduct over the $i$th Rost operations.

We will be mostly interested in one special case: when $j=-k$, $\rho_i^W$ is an 
operation $\textrm{Ch}_k(\bu{W})\to
\textrm{Ch}_{k-i}(X)$ and $\rho^W$ is an operation $\textrm{Ch}(\bu{W})\to 
\textrm{Ch}(X)$.
\end{defn}

The next proposition will require the following two results.

\begin{lemma}\label{fsproper} Let $f:X\to Y$ and $g:Y\to W$ be morphisms 
of (quasi-projective) schemes.  If $g\circ f$ is proper and $f$ is finite and 
surjective, then $g$ is proper.
\end{lemma}
\begin{proof} Take any quasi-projective scheme $T$, and let $f':X\times 
T\to Y\times T$ and $g':Y\times T\to W\times T$ be the induced morphisms.  
Take $Z\subset Y\times T$ a closed subvariety.  Then 
$C=(f')^{-1}(Z)\subset X\times T$ is closed and so $gf(C)\subset 
W\times T$ is closed, as $g\circ f$ is proper.  But $gf(C)=g(Z)$, so 
$g$ is proper, the desired result.
\end{proof}

\begin{cor}\label{XGYGproper}Let $f:X\to Y$ be a proper $G$-equivariant 
morphism of $G$-schemes.  Then the induced morphism $f':X/G\to Y/G$ is 
also proper.
\end{cor}
\begin{proof}Consider the commutative square
\[
\xymatrix{
X \ar[r]^f \ar[d]^{\phi} & Y \ar[d]^{\psi} \\
X/G \ar[r]^{f'} & Y/G
}
\]
Now $\psi$ is finite, hence proper, and so $\psi\circ f=f'\circ\phi$ is 
proper.  The result follows from Lemma \ref{fsproper}.
\end{proof}

\begin{prop}\label{rostfunctoriality}
Let \(f:S\to W\) be a closed embedding of $G$-schemes, and $\bar{f}:\bu{S}\to\bu{W}$ 
the induced morphism.  Let $X=S^G$ and $Y=W^G$.  Then, for any $k,j\in\ZZ ,i\in\NN$, 
the following diagram commutes:
\[
\xymatrix{ 
A_{k,j}(\bu{S}) \ar[r]^-{\rho^{S}_i} \ar[d]_-{\bar{f}_*} 
& A_{k-i,j+i}(X) \ar[d]^-{f_*} \\
A_{k,j}(\bu{W}) \ar[r]^-{\rho^{W}_i} & A_{k-i,j+i}(Y)
}
\]
\end{prop}

\begin{proof}
We have already seen that $\bar{f}$ is proper (Lemma \ref{barf}).  
In light of Proposition \ref{alphafunctoriality} and \cite{merkurjev}, Proposition 53.3, 
it remains to show that $\partial\circ\bar{f}_*=f_*\circ\partial$.  
We consider the triples involved:
\[
\xymatrix{
\bu{S} \ar[r] \ar[d]^{\bar{f}} 
& S/G \ar[d]^{f'} 
& X \ar[l] \ar[d]^f \\
\bu{W} \ar[r] & W/G & Y \ar[l]
}
\]
By Corollary \ref{XGYGproper}, we see that $f'$ is proper; hence we can apply 
\cite{merkurjev}, Proposition 49.33 and we are done.
\end{proof}

\subsection{Pull-backs and the Rost operation}

In general, the Rost operation does not commute with pull-backs.  However, 
we do have such a result in a special case.

\begin{prop}\label{smoothpullbacks} Suppose we have a commutative diagram
\[
\xymatrix{
W \ar[r]^f & Y \\
W^G \ar[u] \ar[r]^g & Y^G \ar[u]
}
\]
where $f$ and $g$ are $G$-equivariant regular closed embeddings with normal 
bundles $V$ and $E$, respectively.  Suppose further that $V^G=E$.  Then this 
gives rise to a commutative diagram
\[
\xymatrix{
\textnormal{Ch}(\bu{Y}) \ar[r]^{\bar{f}^*} \ar[d]^{\rho} & \textnormal{Ch}(\bu{W}) 
\ar[d]^{\rho} \\
\textnormal{Ch}(Y^G) \ar[r]^{g^*} & \textnormal{Ch}(W^G)
}
\]
\end{prop}
\begin{proof}
Recall that by Lemma \ref{U_Wcone}, $\bu{V}$ is the normal bundle of $\bu{W}$ 
in $\bu{Y}$.  Let $\bar{f}:\bu{W}\to\bu{Y}$ be the induced morphism, 
$q_1:\bu{V}\to\bu{W}$ and $q_2:E\to W^G$ the projections.  Then we want 
to prove commutativity in the following diagram:
\[
\xymatrix{
\textnormal{Ch}(\bu{Y}) \ar[r]^{\sigma_{\bar{f}}} \ar[d]^{\rho} \ar@{}[dr]|-{(a)}
& \textnormal{Ch}(\bu{V}) \ar[r]^{(q_1^*)^{-1}} \ar[d]^{\rho} \ar@{}[dr]|-{(b)}
& \textnormal{Ch}(\bu{W}) \ar[d]^{\rho} \\
\textnormal{Ch}(Y^G)\ar[r]^{\sigma_g}
& \textnormal{Ch}(E) \ar[r]^{(q_2^*)^{-1}} & \textnormal{Ch}(W^G)
}
\]

Let $q:\Gm\times\bu{Y}\to\bu{Y}$ and $r:\Gm\times Y^G\to Y^G$ be the projections.  
Then (a) becomes (after rotating by ninety degrees):
\[
\xymatrix{
A(\bu{Y}) \ar[r]^-{\alpha} \ar[d]^-{q^*} \ar@{}[dr]|-{(1)}&
A(\bu{Y}) \ar[r]^-{e} \ar[d]^-{q^*} \ar@{}[dr]|-{(2)}&
A(\bu{Y}) \ar[r]^-{\partial}  \ar[d]^-{q^*} \ar@{}[dr]|-{(3)}&
A(Y^G) \ar[d]^-{r^*} \\
A(\Gm\times\bu{Y}) \ar[r]^-{\alpha} \ar[d]^-{\{ t\}} \ar@{}[dr]|-{(4)}&
A(\Gm\times\bu{Y}) \ar[r]^-{e} \ar[d]^-{\{ t\}} \ar@{}[dr]|-{(5)}&
A(\Gm\times\bu{Y}) \ar[r]^-{\partial}  \ar[d]^-{\{ t\}} \ar@{}[dr]|-{(6)}&
A(\Gm\times Y^G) \ar[d]^-{\{ t\}} \\
A(\Gm\times\bu{Y}) \ar[r]^-{\alpha} \ar[d]^-{\partial} \ar@{}[dr]|-{(7)}&
A(\Gm\times\bu{Y}) \ar[r]^-{e} \ar[d]^-{\partial} \ar@{}[dr]|-{(8)}&
A(\Gm\times\bu{Y}) \ar[r]^-{\partial}  \ar[d]^-{\partial} \ar@{}[dr]|-{(9)}&
A(\Gm\times Y^G) \ar[d]^-{\partial} \\
A(\bu{V}) \ar[r]^-{\alpha} &
A(\bu{V}) \ar[r]^-{e} &
A(\bu{V}) \ar[r]^-{\partial} &
A(E)
}
\]

Some of these are easy to handle.  (1) commutes by Proposition \ref{alphafunctoriality}.  (2) commutes by \cite{merkurjev}, Proposition 53.3 (2).  (4) certainly anticommutes, (5) commutes by \cite{merkurjev}, Propositions 49.7 and 49.16 (we view $e$ as the composition of a push-forward and the inverse to a pull-back).  (6) and (7) anticommute by \cite{merkurjev}, Proposition 49.33, and (9) anticommutes by \cite{merkurjev}, Proposition 49.36.

For (3), we notice that this arises from the following triples:
\[
\xymatrix{
\Gm\times\bu{Y}\ar[r]\ar[d] & \Gm\times Y/G\ar[d] & \Gm\times Y^G \ar[l] \ar[d] \\
\bu{Y} \ar[r] & Y/G & Y^G \ar[l]
}
\]

Certainly, all of the vertical morphisms are flat so we can apply \cite{merkurjev},
Proposition 49.33 (2) to conclude that (3) anticommutes.

In (8), the situation is the following:
\[
\xymatrix{
A(\Gm\times\bu{Y})\ar[r]^-{s_*}\ar[d]^{\partial} 
& A(M) \ar[r]^-{(\phi^*)^{-1}} \ar[d]^{\partial}
& A(\Gm\times\bu{Y}) \ar[d]^{\partial} \\
A(\bu{V})\ar[r]^{r_*} & A(N)\ar[r]^{(\psi^*)^{-1}} & A(\bu{V})
}
\]
where
\[
\xymatrix{
M=(\Gm\times U_Y\times\Aff )/G \ar@{}[dr]|-{\textrm{and}} \ar@/_/[d]_{\phi}
& N=(U_V\times\Aff )/G \ar@/_/[d]_{\psi} \\
\Gm\times\bu{Y} \ar@/_/[u]_{s} & \bu{V} \ar@/_/[u]_{r}
}
\]
are line bundles, with the given maps projections and zero sections.  This gives rise to
\[
\xymatrix{
M \ar[r] \ar@/_/[d]_{\phi} & (U_D\times\Aff )/G \ar@/_/[d]_{\pi} & N \ar[l] \ar@/_/[d]_{\psi} \\
\Gm\times\bu{Y} \ar@/_/[u]_{s} \ar[r] & \bu{D} \ar@/_/[u]_z & \bu{V} \ar@/_/[u]_{r} \ar[l]
}
\]
where $D$ is the deformation variety corresponding to the embedding $W\to Y$.  Now $z$ is proper and $\pi$ is flat, so we can apply \cite{merkurjev}, Proposition 49.33 and we prove the commutativity in (8).

Overall, we had five squares commute and four anticommute, so the diagram is commutative.  Hence (a) is done.

For (b), it suffices to show the commutativity of
\[
\xymatrix{
\Ch{\bu{W}} \ar[r]^{q_1^*} \ar[d]^{\rho} & \Ch{\bu{V}} \ar[d]^{\rho} \\
\Ch{W^G} \ar[r]^{q_2^*} & \Ch{E}
}
\]
By Proposition \ref{alphafunctoriality} and \cite{merkurjev}, Proposition 53.3 (2), it remains to show the commutativity of
\[
\xymatrix{
A(\bu{W}) \ar[r]^{q_1^*} \ar[d]^{\partial} & A(\bu{V}) \ar[d]^{\partial} \\
A(W^G) \ar[r]^{q_2^*} & A(E)
}
\]
We expand to see the triples:
\[
\xymatrix{
\bu{V} \ar[d]^{q_1} \ar[r] & V/G \ar[d]^{g} & V^G=E \ar[l] \ar[d]^{q_2} \\
\bu{W} \ar[r] & W/G & W^G \ar[l]
}
\]
$V/G$ is a vector bundle over $W/G$, hence $g$ is flat and we can again apply \cite{merkurjev},
Proposition 49.33 to finish the proof.
\end{proof}

\begin{cor}\label{transversediagram} Let $W\to X$ be a vector bundle with $W^G=X$.
Let $f:Y\to X$ be a regular closed embedding.  Consider the transverse diagram
\[
\xymatrix{
W|_Y \ar[r] \ar[d] & W \ar[d] \\
Y \ar[r]^f & X
}
\]
Then $\rho^W\circ f^*=\bar{f}^*\circ\rho^{W|_Y}$.
\end{cor}
\begin{proof}
The normal bundle of $W|_Y$ in $W$ is the pull-back of the normal bundle of $Y$ in $X$,
so this diagram satisfies the conditions of Proposition \ref{smoothpullbacks}.
\end{proof}

\subsection{Reduction to the case of cones}

\begin{defn} Let $X=\spec{R}=\spec{(\Rp )}$ be an affine $G$-scheme.  We say a function $f:X\to\Aff$ is \textit{$G$-invariant} if $f\in R_0$.  This definition is extended to arbitrary $G$-schemes in the obvious way.
\end{defn}

We only need one property of these objects.

\begin{lemma}\label{invariantfunction} Let $X$ be a $G$-scheme, $f:X\to\Aff$ a G-invariant function.  Let $Y=V(f)\subset X$.  Then:

(1) The normal cone of $Y$ in $X$ is $L=Y\times\Aff$.

(2) The normal cone of $Y^G$ in $X^G$ is $E=Y^G\times\Aff$.

(3) $L^G=E$.
\end{lemma}
\begin{proof}
Both (1) and (2) are clear.  For (3), we may assume that all schemes are affine, so that $X=\spec{R}=\spec{\Rp}$, $Y=\spec{R/K}$ where $K=fR\cong R$, $X^G=\spec{R/J}$, where $J=I\oplus R_1\oplus\cdots\oplus R_{p-1}$ as usual, and $Y^G=\spec{R/(J+K)}$.

Now 
\[L=\spec{R/K\oplus K/K^2\oplus\cdots}=\spec{R/K[\bar{f}]}\]
where $\bar{f}$ is the image of $f$ under $K\cong R$.  Furthermore, 
\[E=\spec{R/(J+K)\oplus (J+K)/(J+K)^2\oplus\cdots}=
\spec{R/(J+K)[\tilde{f}]}\]
where $\tilde{f}$ is the image of $f$ under $J+K\cong R$.  Finally, since $f\in R_0$, $G$ acts trivially on $f$, meaning that
\(L^G=\spec{R/(J+K)[\tilde{f}]}\),
the desired result.
\end{proof}

\begin{cor}\label{reductionvectorbundles} Let $W$ be a $G$-scheme, 
$X=W^G$.  Let $C$ be the normal cone of $X$ in $W$.  Then
\[\rho^W([\bu{W}])=\rho^C([\bu{C}])\in\Ch{X}\]
\end{cor}
\begin{proof}
Let $D$ be the deformation variety.  The closed embeddings $f:C\to D$ and $g:W\to D$ are both regular.  Consider the following diagram:
\[
\xymatrix{
\Ch{\bu{C}} \ar[d]^{\rho^C} & \Ch{\bu{D}} \ar[r]^{\bar{g}^*} \ar[d]^{\rho^D} \ar[l]_{\bar{f}^*} 
& \Ch{\bu{W}} \ar[d]^{\rho^W} \\
\Ch{X} & \Ch{X\times\Aff} \ar[l] \ar[r] & \Ch{X}
}
\]
Note that $C\subset D$ and $W\subset D$ are both defined by $G$-invariant
functions.  Therefore, by Proposition \ref{smoothpullbacks} and Lemma \ref{invariantfunction}, both small squares commute.  Consider $[\bu{D}]\in\Ch{\bu{D}}$.  If we follow it left and down, it is mapped to $\rho^C([\bu{C}])$.  If we follow it right and down, it is mapped to 
$\rho^W([\bu{W}])$.  But these two must be equal because the bottom row is 
equality.
\end{proof}

\section{The map $P_X$ and the Steenrod operations}\label{steenroddef}

In this section we finally describe the construction of the Steenrod operations.  The last component we need is the map $P_X$, which can be viewed as a generalization of the map $P^n_G$ in \cite{brosnan}.

Let us set up some notation.  Take a scheme $X$ with trivial $G$-action; consider $R^p(X)$ and the associated $G$-action as before.  We have seen that $X=(R^p(X))^G$; we set $W=R^p(X)$ and use this example from now on.  We simplify notation in what follows: instead of $U_{R^p(X)}$ and $\bu{R^p(X)}$, whenever convenient, we use the notation $U$ and $\bu{}$.  We also write $\rho$ and $\rho_{i}$ to stand for $\rho^{R^p(X)}$ and $\rho^{R^p(X)}_{i}$.  If $X$ is not clear from context, we will abuse notation and denote these by $U_X,\bu{X},\rho^X$ and $\rho^X_i$.

\subsection{Twists of a $G$-action}

\begin{defn}Let $W=\spec{R}$ be an affine scheme with a $G$-action, so that $R=\Rp$.
Suppose $R_i\cong R_j$ for all $i,j=1,2,\ldots ,p-1$.  Then for each $k=1,2,\ldots, 
p-1$ we have an isomorphism of $R$:
\[h_k:(r_0,r_1,\ldots ,r_{p-1})\mapsto (r_0,r_{k'},\ldots ,r_{k'(p-1)})\]
where $k'$ is the inverse of $k$ (mod $p$).  This induces an isomorphism 
$\eta_k:W\to W$, known as a \textit{$k$-twist of $W$.}  

If in addition $W$ has a canonical $G$-action, and we act on $W$ by $G$ and 
then apply a $k$-twist, we say that we have \textit{twisted the action $k$ times} 
and that $G$ acts on $W$ with \textit{weight $k$}.  Schemes on which we can apply 
a $k$-twist include $X^p$ and $R^p(X)$ for any scheme $X$, along with cones and
vector bundles.  The last two have a canonical $G$-action, so we can talk about
twisting the action and the weight of the action.  Note that this corresponds with
the standard definition of twisting a $G$-action, which is composing the action map
$G\times W\to W$ with $G\times W\to G\times W$ given by $(g,x)\mapsto (g^k,x)$.
\end{defn}

\begin{lemma}\label{alphaeta}Let $X$ be a scheme, and let $L=(U_X\times\Aff )/G$
be the canonical line bundle over \bu{X}.  Then for any $i,j\in\ZZ$ and 
$k=1,2,\ldots ,p-1$,
\begin{enumerate}

\item\(k\cdot (\eta_k)_*\circ\alpha_{\bu{}}=\alpha_{\bu{}}\circ (\eta_k)_*:
A_{i,j}(\bu{})\to A_{i,j+1}(\bu{})\)

\item\(k\cdot (\eta_k)_*\circ\ e(L)=e(L)\circ (\eta_k)_*:
A_{i,j}(\bu{})\to A_{i-1,j+1}(\bu{})\)

\item\(\partial^{\bu{X}}_X=\partial^{\bu{X}}_X\circ (\eta_k)_*:
A_{i,j}(\bu{})\to A_{i-1,j}(X)\)
\end{enumerate}

\end{lemma}
\begin{proof}(1) Consider the diagram
\[
\xymatrix{
C_m(\bu{}) \ar[d]^{(\eta_k)_*} \ar[r]^{\alpha_{\bu{}}} & C_m(\bu{}) 
\ar[d]^{(\eta_k)_*} \\
C_m(\bu{}) \ar[r]^{\alpha_{\bu{}}} & C_m(\bu{})
}
\]
The $(x,y)$-component in either direction is nonzero only if $y=\eta_k(x)$.  
We can write $F(x)=\Rp$ and $F(y)=S_0\oplus S_1\oplus\cdots\oplus S_{p-1}$ where
$S_i=R_{ik}$ for each $i$.  Now (in the notation of Example \ref{tp-r}) since $R_i=b_x^iR_0$, we have that $S_i=b_x^{ik}R_0$, 
and so $a_y=(b_x^k)^p=a_x^k$.  The result now follows.

(2) $(\eta_k)^*(L)=L^{\otimes k}$, hence \((\eta_k)_*\circ\ e(L^{\otimes k})=e(L)
\circ (\eta_k)_*\) (by \cite{merkurjev}, Proposition 53.3), giving the result.

(3) follows from the commutativity of
\[
\xymatrix{
\bu{X} \ar[r] \ar[d]^{(\eta_k)_*} & R^p(X)/G \ar[d]^{(\eta_k)_*} & X
\ar@{=}[d] \ar[l] \\
\bu{X} \ar[r] & R^p(X)/G & X \ar[l]
}
\]
and \cite{merkurjev}, Proposition 49.33 (1).
\end{proof}

\begin{prop}$\rho_i^X:\bu{X}\to X$ is trivial unless $p-1|\,i$.
\end{prop}
\begin{proof}Consider the diagram (subscripts omitted for brevity):
\[
\xymatrix{
\Ch{\bu{X}} \ar[r]^{\alpha_{\bu{X}}} \ar[d]^{(\eta_k)_*} &
A(\bu{X}) \ar[r]^{e(L)^{i-1}} \ar[d]^{(\eta_k)_*} &
A(\bu{X}) \ar[r]^{\partial^{\bu{X}}_X} \ar[d]^{(\eta_k)_*} &
\Ch{X} \ar@{=}[d] \\
\Ch{\bu{X}} \ar[r]^{\alpha_{\bu{X}}} &
A(\bu{X}) \ar[r]^{e(L)^{i-1}} &
A(\bu{X}) \ar[r]^{\partial^{\bu{X}}_X} &
\Ch{X}
}
\] 
Both the top and bottom compositions equal $\rho_i$, and hence by Lemma 
\ref{alphaeta}, $k^i\cdot\rho_i=\rho_i$.  Since this is true for any 
$k\in\Zps$, it must be that $\rho_i=0$ unless $p-1|\,i$.
\end{proof}

\subsection{The map $P_X$}

For any scheme $X$, we set $Z(X;\PP^1)$ to be the subgroup of
$Z(X\times\PP^1)$ generated by the classes of closed subvarieties in
$X\times\PP^1$ that are dominant over $\PP^1$.  For a rational point 
$a\in\PP^1$ and a variety $V\subset X\times\PP^1$ mapping dominantly to 
$\PP^1$, let $V(a)$ denote the fiber above $a$.  We associate to this the 
cycle $[V(a)]$.  This definition is extended to $Z(X;\PP^1)$ by 
linearity.

\begin{prop}\label{ratequivcycles}
Let $\gamma_0$ and $\gamma_{\infty}$ be two cycles on a scheme $X$.  Then 
the classes of $\gamma_0$ and $\gamma_{\infty}$ in $\CH{X}$ are equal if 
and only if there exists a cycle $\delta\in Z(X;\PP^1)$ such that 
$\gamma_0=\delta (0)$ and $\gamma_{\infty}=\delta (\infty )$
\end{prop}

\begin{proof}
See \cite{brosnan}, Proposition 4.2  
\end{proof}

Take a closed subvariety $Z\subset X$.  Then $R^p(Z)\subset R^p(X)$ is a
closed subvariety such that $[R^p(Z)]\subset Z(R^p(X))^G$.  We may
extend this by linearity to any cycle $\gamma\in Z(X)$ to get
a $G$-invariant \(R^p(\gamma )\in Z(R^p(X))\).

By Corollary \ref{torsorisomorphism} we get a cycle $\bar{\gamma}\in
Z(\bu{})$ such that $\pi^*(\bar{\gamma})=R^p(\gamma )|_{U}$.  We hence
have a map $Z(X)\to Z(\bu{})$, taking $\gamma\mapsto\bar{\gamma}$.

\begin{lemma}
If $\gamma_0$ and $\gamma_{\infty}$ are rationally equivalent cycles in
$Z(X)$, then $\bar{\gamma}_0$ and $\bar{\gamma}_{\infty}$ are rationally
equivalent cycles in $Z(\bu{})$.
\end{lemma}

\begin{proof}
Let $Z\subset X\times \PP^1$ be a closed subvariety dominant over $\PP^1$.
Consider the following fiber diagram:
\[
\xymatrix@=12pt{
&& \wt{Z} \ar@{}|-{\bigcap}[d] \ar@{}|-{\subset}[r] & R^p(X)\times\PP^1
\ar@{}|-{\bigcap}[d] \\
&& R^p(Z) \ar@{}|-{\subset}[r] & R^p(X\times\PP^1) 
\ar@{}|-{=}[r] & R^p(X)\times R^p(\PP^1) \\
}
\]
where $\widetilde{Z}=(R^p(X)\times\PP^1)\cap R^p(Z)\subset
R^p(X)\times\PP^1$ is the fiber product. Notice that over an algebraic closure,
\(\widetilde{Z}=Z\times_{\PP^1}Z\times_{\PP^1}\cdots\times_{\PP^1}Z\).  By
descent theory, we see that
$\widetilde{Z}\to\PP^1$ is flat.  Therefore, every irreducible
component of $\widetilde{Z}$ is dominant over $\PP^1$, i.e.\ the cycle
$[\widetilde{Z}]$ belongs to $Z(R^p(X);\PP^1)$.  By linearity, this 
construction extends to a map
\[Z(X;\PP^1)\to 
Z(R^p(X);\PP^1),\,\,\,\,\,\delta\mapsto\widetilde{\delta}\]

Let $G$ act on $R^p(X)\times\PP^1$ by the standard action on $R^p(X)$ and trivially on $\PP^1$.  By Proposition \ref{ratequivcycles}, we have a cycle $\delta\in Z(X;\PP^1)$ such that $\delta (0)=\gamma_0$ and $\delta(\infty )=\gamma_{\infty}$.  The associated cycle $\widetilde{\delta}\in Z(R^p(X);\PP^1)$ is $G$-invariant. Since $U\times\PP^1$ is a $G$-torsor over $\bu{}\times\PP^1$, the restriction of the cycle $\widetilde{\delta}$ on $U\times\PP^1$ gives rise to a well-defined cycle $\widehat{\delta}\in Z(\bu{};\PP^1)$, such that
\begin{equation}
\label{part1}
\wt{\delta}|_{U\times\PP^1}=q^*(\widehat{\delta})
\end{equation}
where
\(q:U\times\PP^1\to\bu{}\times\PP^1\) is the canonical morphism.

Let \(Z\subset \bu{}\times\PP^1\) be a closed subvariety dominant over
$\PP^1$.  We have 
\(\pi^{-1}(Z(a))=q^{-1}(Z)(a)\)
for any rational point $a$ of $\PP^1$, where $\pi :U\to \bu{}$ is the
canonical morphism.  It follows from \cite{merkurjev}, Proposition 57.7 that 
\begin{equation}
\label{9.6}
\pi^*(\eta (a))=(q^*\eta )(a)
\end{equation}
for every cycle $\eta\in Z(\bu{};\PP^1)$.

Let $\delta =\sum n_i[W_i]$.  Then applying (\ref{9.6}) for $\eta
=\widehat{\delta}$, we get from (\ref{part1}) that
\bears
\pi^*(\widehat{\delta}(a)) &=&
(q^*\widehat{\delta})(a)=\widetilde{\delta}|_{U\times\PP^1}(a) \\
&=& \sum n_i[\wt{W_i}]|_{U\times\PP^1}(a) \\
&=& \sum n_i[R^p(W_i(a))]|_U \\
&=& R^p(\delta (a))|_U \\
&=& \pi^*(\overline{\delta (a)})
\eears
Since $\pi^*$ is injective, $\widehat{\delta}(a)=\overline{\delta (a)}$ in
$Z(\bu{})$.  In particular,  $\widehat{\delta}(0)=\overline{\delta
(0)}=\bar{\gamma}_0$ and $\widehat{\delta}(\infty )=\overline{\delta (\infty
)}=\bar{\gamma}_{\infty}$, i.e.\ the cycles $\bar{\gamma}_0$ and
$\bar{\gamma}_{\infty}$ are rationally equivalent, as desired.
\end{proof}
By this lemma, we have a well-defined map
\[
P_X:\textrm{CH}(X)\to \textrm{CH}(\bu{}),\,\,\,\,\,[\gamma ] \mapsto 
[\bar{\gamma}]
\]
If we like, we can also consider $P_X$ as a map $\textrm{CH}_n(X)\to \textrm{
CH}_{np}(\bu{})$ for any $n$.  When $X$ is clear from context, we 
will simplify notation and denote this map by $P$.  Note that this map 
is not a homomorphism.  

\begin{prop}\label{Pfunctoriality}
Let $f:X\to Y$ be a closed embedding.  Then \(\bar{f}_*\circ P_X=P_Y\circ f_*\).  Moreover, if $f$ is a regular embedding, \(\bar{f}^*\circ P_Y=P_X\circ\bar{f}^*\).
\end{prop}

\begin{proof} 
Note that this proposition makes sense in view of Lemma \ref{barf}.  The first stated commutativity is easy, as it holds on the level of cycles.  For the second, we have by Lemma \ref{U_Wcone} that the normal bundle of $\bu{X}$ in $\bu{Y}$ is $\bu{N}$, where $N$ is the normal bundle of $X$ in $Y$.  Now the map $f^*$ is defined as the composition $(\pi^*)^{-1}\circ\sigma$ where $\sigma:\CH{\bu{Y}}\to\CH{\bu{N}}$ is the deformation homomorphism and $\pi:\bu{N}\to\bu{Y}$ is the projection which induces an isomorphism on the Chow groups.

The map $\sigma$ is defined on cycles by $\sigma([Z])=[C_{Z\cap X}Z]$, where $C_{Z\cap X}Z$ is the normal cone of $Z\cap X$ in $Z$.  Hence $\sigma$ commutes with $P$ on the level of cycles.  Furthermore, it is easy to see that $P$ commutes with $\pi^*$ on the level of cycles, and so we are done.
\end{proof}

\noindent \textbf{Note.}  $P$ may not commute with flat pull-backs in general as $\bar{f}$ may not even be defined for a flat morphism $f$.

\subsection{Definition of the Steenrod operations}

\begin{defn} We define the \textit{k-th Steenrod operation mod p} to be 
\[S_k^X=(-1)^{n+k}\rho_{(p-1)(n+k)}\circ P:\textrm{Ch}_n(X)\to \textrm{
Ch}_{n-(p-1)k}(X)\]
\end{defn}
By taking the coproduct over all of these, we can define the \textit{total 
Steenrod operation} $S^X:\Ch{X}\to\Ch{X}$.

\begin{prop} $S^X$ is a homomorphism.
\end{prop}
\begin{proof}
By a restriction-corestriction argument, we may assume that $F$ contains 
all $p$-th roots of unity.  Let $\pi :U\to\bu{}$ be the projection.  
Choose a distinguished $p$-th root of unity $\xi$ of $F$ and let $\sigma$ 
be the automorphism of $Z(U)$ associated to $\xi$.

For any two cycles $\gamma =\sum n_i[Z_i]$ and $\delta =\sum m_i[Z_i]$ on 
$X$ we have
\[
\pi^*(\overline{\gamma 
+\delta})-\pi^*(\overline{\gamma})-\pi^*(\overline{\delta})
=\sum n_im_j(1+\sigma+\cdots +\sigma^{p-1})[Z_i\times 
Z_j]|_{U}\in Z(U)
\]
Since $\pi^*\circ\pi_*=1+\sigma+\cdots +\sigma^{p-1}$ by Proposition 
\ref{gtorsorprops} (3), we have
\[
P(\gamma +\delta )-P(\gamma )-P(\delta )\in \textrm{Im}(\pi_*)\subset\CH{U}
\]
so it suffices to show $\rho\circ\pi_*=0$.  Since $\rho 
=\partial\circ e\circ\alpha$, it is certainly enough to show 
$\alpha\circ\pi_*=0$.

Consider the fiber diagram
\[
\xymatrix{
G\times U \ar[r] \ar[d] & U \ar[d]^{\pi} \\
U \ar[r]^{\pi} & \bu{}
}
\]
where \(G\times U\) is the fiber product, as $U\to\bu{}$ is a $G$-torsor.  
By Proposition \ref{alphanaturality} (2), we get the following commutative 
diagram:
\[
\xymatrix{
\Ch{U} \ar[r]^{\alpha_{G\times U}} \ar[d]^{\pi_*} 
	& \Ch{U} \ar[d]^{\pi_*} \\
\Ch{\bu{}} \ar[r]^{\alpha_U} & \Ch{\bu{}}
}
\]

But by Proposition \ref{alphanaturality} (1), $\alpha_{G\times U}=0$.
Hence $\alpha\circ\pi_*=0$ and we are done.
\end{proof}

\begin{prop}\label{steenrodclosedembedding}
Let $i:X\to Y$ be a closed embedding.  Then $i_*\circ S^Y=S^X\circ i_*$.
\end{prop}
\begin{proof}
This is immediate from Propositions \ref{rostfunctoriality} and 
\ref{Pfunctoriality}.
\end{proof}

\section{Computation of the Steenrod operations}\label{compsteenrod}

\subsection{Some generalizations of Chern classes}

\begin{defn}
Let $V$ be a vector bundle of rank $n$ on $X$ and let
\(1+c_1+c_2+\cdots +c_n=\prod_{i=1}^{n}(1+x_i)\)
so that the $x_i$ are the Chern roots of $V$.

\begin{enumerate}
\item We set 
\(b(V)=\prod_{i=1}^n(1+x_i^{p-1})=1+b_1+b_2+\cdots\)
and call it the \textit{b-class}, or \textit{multiplicative Chern class} of $V$.

\item Set 
\(\omega(V)=\prod_{i=1}^n(1-x_i^{p-1})=1+\omega_1+\omega_2+\cdots\)
and call it the \textit{omega-class} of $V$.

\end{enumerate}

\end{defn}

\noindent \textbf{Remarks:}
\begin{enumerate}
\item The operation $b$ is multiplicative, i.e.\ $b(V\oplus W)=b(V)b(W)$ for any vector bundles $V$ and $W$.  Furthermore, $b(L)=1+c_1(L)^{p-1}$ for a line bundle $L$.  By the splitting principle, these two properties determine $b$.
\item If $p=2$, $\omega(V)=1+c_1+c_2+\cdots =c(V)=b(V).$
\item Since the series $b$ lies in the subring of $R$ of elements of 
degree divisible by $p-1$, $b$ is of the form
\(b=1+b_{p-1}+b_{2(p-1)}+\cdots \)
with $b_i$ of degree $i$.  The same can be said of $\omega$.  In fact, it is easy to see that $\omega_{k(p-1)}=(-1)^kb_{k(p-1)}.$
\end{enumerate}

\begin{defn} Let $M\to X$ be a $G$-vector bundle with $M^G=X$.  Suppose $M$
has a filtration by $G$-subbundles
\[0=M_0\subset M_1\subset\cdots\subset M_n=M\]
with each quotient $L_i=M_i/M_{i-1}$ a line bundle.  For $i=1,\ldots n$, let 
$r_i$ be the weight of $L_i$ as a $G$-module.  Then we set
\[\mu (M)=\prod_{i=1}^n(r_i+c_1(L_i))\]
and call it the \textit{mu-class of $M$}.
\end{defn}

\begin{lemma}\label{mu-rho} Let $j:V\to E$ be a closed embedding of $G$-vector bundles over a scheme $Y$, with $V^G=E^G=Y$.  Suppose $E/V$ has a filtration
as in the above definition; finally, let $i:\bu{V}\to\bu{E}$ be the induced closed
embedding.  Then
\[\rho^V\circ i^*=\mu (E/V)\circ\rho^E\]
\end{lemma}
\begin{proof}
Assume that $V\subset E$ is a subbundle of corank 1; the general result will follow 
by induction.  Take $\gamma\in\Ch{\bu{E}}$.  Then by the projection formula and
Proposition \ref{rostfunctoriality}, 
\begin{equation}\label{mu-rho-eqn}
\rho^Vi^*(\gamma )=\rho^Ei_*i^*(\gamma )
=\rho^Ei_*(1\cdot i^*(\gamma ))=\rho^E(i_*(1)\cdot\gamma )
=\rho^E(c_1(L)(\gamma))
\end{equation}
where $L$ is the line normal bundle of $\bu{V}\subset\bu{E}$.

Now consider the normal bundle $E/V$ of $V$ in $E$; let $r$ be its weight as a 
$G$-module.  We may write it as $E/V\otimes\Aff$ where $G$ acts trivially on 
$E/V$ and with weight $r$ on \Aff.  Then we have that $L=E/V\otimes L_E^{\otimes 
r}$, where $L_E=(U_E\times\Aff )/G$ is the canonical line bundle over $\bu{E}$.

Hence, continuing from (\ref{mu-rho-eqn}), 
\bears
\rho^V\circ i^*(\gamma )=\rho^E(c_1(L)(\gamma))
&=& \rho^E(rc_1(L_E)(\gamma )+c_1(E/V)(\gamma )) \\
&=& r\rho^E(\gamma)+c_1(E/V)\rho^E(\gamma) \\
&=& (r+c_1(E/V))\rho^E(\gamma)
\eears
completing the proof.
\end{proof}

\subsection{The main computation}

At this point we need to introduce the concept of equivariant Chow groups.  We will not be using them much, so we will simply refer the interested reader to \cite{edidin} or \cite{brosnan} for further reading.  Two well-known facts from this theory will be of most importance:

\noindent \textbf{Fact 1.} Set $l=c_1(L)\in A^1B\mu_p$, where $L$ is the canonical equivariant line bundle over $\textnormal{pt}_F$.  Then if $X$ has trivial $G$-action, $\ChG(X)=\Ch{X}[l]$ (cf.\ \cite{brosnan}, Theorem 7.1.)

\noindent \textbf{Fact 2.} If $T\to T/G$ is a $G$-torsor, then $\ChG(T)=\Ch{T/G}$.

\begin{defn} Let $E\to X$ be a $G$-vector bundle on a variety $X$ with trivial $G$-action such that $E^G=X$.  Let $j:E\setminus X\to E$ and $s:E\to X$ be the inclusion and projection, respectively.  We set $\Phi=\Phi_{E,X}:\ChG(X)\to\Ch{\bu{E}}$ to be the composition
\[\Phi:\ChG(X)\stackrel{s^*}{\to}\ChG(E)\stackrel{j^*}{\to}\ChG(E\setminus X)=\Ch{\bu{E}}\]
\end{defn}

\noindent\textbf{Note.} Under the association $\ChG(X)=\Ch{X}[l]$ from Fact 1, 
\[\Phi(\delta\cdot l^i)=c_1(L_E)^i(\psi^*(\delta))\]
where $L_E$ is the canonical line bundle on $\bu{E}$ and $\psi:\bu{E}\to X$ is the projection.

\noindent\textbf{Notation.} For any element $\sigma=a_0+a_1l+\cdots +a_nl^n\in\ChG(X)$, we write $\epsilon(\sigma)=a_0+a_1+\cdots +a_n\in\Ch{X}$.

The following Proposition is the main computation.

\begin{prop}\label{calc} Let $q:V\to X$ be a $G$-vector bundle over a variety $X$ with $V^G=X$.  Take a class $\sigma\in\ChG(X)$.  Then under the composition
\[\ChG(X)\stackrel{\Phi}{\to}\Ch{\bu{V}}\stackrel{\rho}{\to}\Ch{X}\]
$\sigma$ is mapped to $\mu(V)^{-1}(\epsilon(\sigma))$.
\end{prop}
\begin{proof}
We may assume $\sigma=\eta\cdot l^a$ where $\eta\in\Chk{k-1}{X}$ for some $k>0$ and $a\in\NN\cup \{ 0\}$.  By the note above, it remains to compute $\rho(c_1(L_V)^a(\psi^*\eta))$ where $L_V$ is the canonical line bundle on $\bu{V}$ and $\psi:\bu{V}\to X$ is the projection.  We proceed by induction on rank $V$.  

Suppose first that rank $V=1$; then $V=X\times\Aff$ where $G$ acts on $\Aff$ with some weight, say $r$.  We have that $V\setminus X=X\times\Gm$ and $\bu{V}=(V\setminus X)/G=(X\times\Gm)/G=X\times\Gm$.  Furthermore, $L_V$ comes from $(\Gm\times\Aff)/G\to\Gm/G$, which is a trivial line bundle.  Hence $L_V$ itself is trivial and so $\rho(c_1(L_V)^a(\psi^*\eta))=0 $ unless $a=0$, in which case it equals $\psi^*\eta$.  We now compute each step of the calculation of $\rho(\psi^*\eta)$.

1. \(\alpha:\Chk{k}{X\times\Gm} \to A_{k,-k+1}(X\times\Gm)\)

From the right exact sequence
\[\Chk{k}{X}\stackrel{0}{\to}\Chk{k}{X\times\Aff}\to\Chk{k}{X\times\Gm}\to 0\]
we find that $\Chk{k}{X\times\Gm}\cong\Chk{k}{X\times\Aff}\cong\Chk{k-1}{X}$.  Note that under this identification, $\psi^*\eta\in\Chk{k}{X\times\Gm}$ corresponds to $\eta\in\Chk{k-1}{X}$.

Similarly, from the split short exact sequence
\[0\to A_{k,-k+1}(X\times\Aff)\to A_{k,-k+1}(X\times\Gm)\stackrel{\partial}{\to}\Chk{k-1}{X}\to 0\]
we get that $A_{k,-k+1}(X\times\Gm)\cong A_{k,-k+1}(X\times\Aff)\oplus\Chk{k-1}{X}$.  We consider $\alpha:\Chk{k-1}{X}\to A_{k,-k+1}(X\times\Aff)\oplus\Chk{k-1}{X}$ and want to find the image of the class $\eta$.

Let $R=F[t, t^{-1}]$.  The canonical $G$-action on $\Gm$ corresponds to the canonical grading $R=\Rp$.  In the twisted grading $R=R'_0\oplus R'_1\oplus\cdots\oplus R'_{p-1}$, $R_i=R'_{ir}$.  Hence $R'_{1}=R_{r'}$ where $r'$ is the inverse of $r$ mod $p$.  Thus $\alpha$ is represented by $t^{r'}$ where $t$ is the standard coordinate on $\Aff$, and so $\alpha(\eta)=(t^{r'},\eta).$

2. $e(L_V)^{j-1}:A_{k,-k+1}(X\times\Gm)\to A_{k-j+1,-k+j}(X\times\Gm)$

We have already noted that $L_V$ is a trivial line bundle.  Thus $e(L_V)^{j-1}(\alpha(\eta))$ is nonzero only when $j=1$, and in that case equals $(t^{r'},\eta)$.

3. The last step is to compute $\partial(t^{r'},\eta)$ which equals $r'\eta$.  Hence $\rho^V(\psi^*\eta)=r'\eta$.  Finally, we see that $\mu(V)=r+c_1(V)=r$, so we are done in this case.

Now suppose the result is true for vector bundles of rank $n-1$.  Take a vector bundle $V$ of rank $n$ over $X$ and suppose it has a filtration $0=V_0\subset V_1\subset\cdots\subset V_{n-1}\stackrel{i}{\hookrightarrow} V_n=V$, with each quotient $L_i=V_i/V_{i-1}$ a line bundle on which $G$ acts with weight $r_i$.  Let $\bar{\iota}^*:\bu{V_{n-1}}\to\bu{V_n}$ and $\psi':\bu{V_{n-1}}\to X$ be the obvious morphisms.  Note that 
\begin{equation}\label{iotaeqn}
\bar{\iota}^*(c_1(L_V)^{a}(\psi^*\eta))=c_1(L_{n-1})^{a}(\psi'^*\eta)
\end{equation}
and by Lemma \ref{mu-rho},
\begin{equation}\label{indeqn}
\rho^{V_{n-1}}(\bar{\iota}^*(c_1(L_V)^{a}(\psi^*\eta)))=\mu(V/V_{n-1})\rho^V(c_1(L_V)^{a}(\psi^*\eta))
\end{equation}
The right-hand side of (\ref{indeqn}) equals $(r_n+c_1(L_n))\rho^V(c_1(L_V)^{a}(\psi^*\eta))$.  By our inductive hypothesis and (\ref{iotaeqn}), the left-hand side equals $(\prod_{i=1}^{n-1}(r_i+c_1(L_i)))^{-1}\eta$.  Multiplying both sides by $(r_n+c_1(L_n))^{-1}$ now gives the result.
\end{proof}

\noindent\textbf{Notation.} Let $Z\subset X$ be a subvariety, $C\to Z$ be a cone and $s:E\to X$ a vector bundle such that $C\subset E$.  Consider the class $[C]\in\ChG(E)$; set $C^G_E=(s^*)^{-1}([C])$ to be the corresponding element in $\ChG(X)$.  Finally, set $C_E\in\Ch{X}$ to equal $\epsilon(C^G_E)$.

\begin{cor}\label{omegacor} Let $X$ be a variety, $X\to W$ a closed embedding of $X$ into a smooth scheme $W$ with $\dim W=e$.  Let $i:Z\to X$ be a closed subvariety and $C$ the normal cone of $Z\to R^p(Z)$.  Let $E$ be the restriction to $X$ of the normal bundle of $W\to R^p(W)$.  Then
\[\rho^X([\bu{Z}])=(-1)^e\omega(-T_W|_X)(C_E)\]
\end{cor}
\begin{proof}First note that by Corollary \ref{reductionvectorbundles} and Proposition \ref{rostfunctoriality},
\[\rho^X([\bu{Z}])=i_*\rho^Z([\bu{Z}])=i_*\rho^C([\bu{C}])=\rho^D([\bu{C}])=\rho^{E}([\bu{C}])\]
where $D$ is the normal cone of $X\to R^p(X)$.  Second, note that $\Phi(C^G_E)=[\bu{C}]\in\Ch{\bu{E}}$, so 
$\rho^{E}([\bu{C}])=\mu(E)^{-1}\epsilon(C^G_E)=\mu(E)^{-1}(C_E)$ by Proposition \ref{calc}.  Therefore it suffices to calculate $\mu(E)$.

Set $H$ to be the cokernel of the embedding $F\to F_p=F[t]/(t^p-1)$; then $E=T_W|_X\otimes H$.  Note that $e$ is the rank of $T_W|_X$ as a vector bundle over $X$.  Suppose $T_W|_X$ has a filtration with quotients $M_i$ ($i=1,2,\ldots e$) and $H$ has a filtration with quotients $H_j$ ($j=1,2,\ldots p-1$); note that the $H_j$ are trivial line bundles.  Hence $T_W|_X\otimes H$ has a filtration with quotients $M_i\otimes H_j$, where $G$ acts with weight $j$, and so
\allowdisplaybreaks{
\bears
\mu(T_W|_X\otimes H) &=& \prod_{j=1}^{p-1}\prod_{i=1}^{e}(j+c_1(M_i\otimes H_j)) \\
&=& \prod_{j=1}^{p-1}\prod_{i=1}^{e}(j+c_1(M_i)+c_1(H_j)) \\
&=& \prod_{i=1}^{e}\prod_{j=1}^{p-1}(j+c_1(M_i)) \\
&=& \prod_{i=1}^{e}(-1+c_1(M_i)^{p-1}) \\
&=& (-1)^{e}\prod_{i=1}^{e}(1-c_1(M_i)^{p-1}) \\
&=& (-1)^e\omega(T_W|_X)
\eears
}
\end{proof}

The following lemma is done in \cite{brosnan} but not explicitly stated in this form.

\begin{lemma}\label{brolemma} Let $Z\subset X$ be a subvariety of dimension $n$ and let $X\to W$ be a closed embedding of $X$ into a smooth scheme $W$.  Let $C$ be the normal cone of $Z\to R^p(Z)$ and $E$ the restriction to $X$ of the normal bundle $W\to R^p(W)$.  Then in the expression 
\[C^G_E=\sum_{i=1}^{r}a_{i}l^{i}\]
all terms with $i$ not divisible by $p-1$ are zero.  In particular, we may write $C_E=\gamma_0+\gamma_1+\cdots +\gamma_m$ where each $\gamma_i\in\Chk{n-i(p-1)}{X}$.
\end{lemma}
\begin{proof}
By a restriction-corestriction argument, we may assume our ground field $F$ contains all roots of unity, so that $R^p(X)=X^p$.
Set $e=\dim W$, let $\gamma=[Z]\in\Chk{n}{X}$ and let $\Delta:W\to W^p$ be the diagonal embedding.  By \cite{brosnan}, Proposition 5.2, there exists a well-defined class $\gamma^p_G\in\ChG_{np}(X^p)$ which comes from $\gamma^{\times p}\in\Chk{np}{X^p}$.  Next, note that 
$\Delta^{!}_{G}(\gamma^p_G)=C^G_E\in\ChG(X)$ by definition (\cite{brosnan}, 3.1.1 and \cite{fulton}, 6.2).  Then \cite{brosnan}, Theorem 8.3 gives the result.

Finally, note that each $a_{i(p-1)}\in\Chk{pn+(p-1)(i-e)}{X}$, so for each $i$ we set 
\[\gamma_i=a_{(e-n-i)(p-1)}\in\Chk{n-i(p-1)}{X}\]
and we get that $C_E=\sum\gamma_i$.
\end{proof}

\begin{cor}\label{corcalc} Let $X$ be a variety, $X\to W$ a closed embedding of $X$ into a smooth scheme $W$ with $\dim W=e$.  Let $Z\subset X$ be a closed subvariety of dimension $n$ and $C$ the normal cone of $Z\to R^p(Z)$.  Let $E$ be the restriction to $X$ of the normal bundle of $W\to R^p(W)$ and write $C_E=\sum\gamma_i$ as in Lemma \ref{brolemma}.  Set $\wt{C}_E=\sum(-1)^{e+n+i}\gamma_i$.  Then $S^X([Z])=b(-T_W|_X)(\wt{C}_E)$.
\end{cor}
\begin{proof}Set $\tilde{\gamma_i}=(-1)^{e+n+i}\gamma_i$ for each $i$.  Then by Corollary \ref{omegacor} and Lemma \ref{brolemma},
\bears
S^X_k([Z]) &=& (-1)^{n+k}\rho^X_{(p-1)(n+k)}P([Z]) \\
&=& (-1)^{n+k}\rho^X_{(p-1)(n+k)}(\bu{Z}) \\
&=& (-1)^{n+k}(-1)^e\sum_{i+j=k}\omega_{i(p-1)}(-T_W|_X)(\gamma_j) \\
&=& (-1)^{n+k}(-1)^e\sum_{i+j=k}(-1)^{i}b_{i(p-1)}(-T_W|_X)(-1)^{e+n+j}(\tilde{\gamma_j}) \\
&=& \sum_{i+j=k}b_{i(p-1)}(-T_W|_X)(\tilde{\gamma_j})\\
&=& (b(-T_W|_X)(\wt{C}_E))_k
\eears
and hence $S^X([Z])=b(-T_W|_X)(\wt{C}_E)$, as desired.
\end{proof}

\subsection{An explicit presentation of the class $C^G_E$}

\begin{prop}\label{subcone} Let $q:E\to X$ be a vector bundle of rank $n$.  Let $i:Z\to X$ be a closed subvariety and $C\to Z$ a cone, such that 
\[
\xymatrix{
C \ar[r] \ar[d] & E \ar[d]^q \\
Z \ar[r]^i & X
}
\]
commutes.  Then
\[C^G_E=\sum_{j+k=n}c^G_j(E)i_*s_k^G(C)\in\ChG(X)\]
\end{prop}
\begin{proof}Consider the commutative diagram
\[
\xymatrix{
\PP(C\oplus\one) \ar[r]^j \ar[d]^l & \PP(E\oplus\one) \ar[d]^r \\
Z \ar[r]^i & X
}
\]
Let $L$ and $L'$ be the canonical line bundles on $\PP(E\oplus\one)$ and $\PP(C\oplus\one)$, so that $L'=i^*L$.  By the equivariant projective bundle theorem, $j_*[\PP(C\oplus\one)]\in\ChG(\PP(E\oplus\one))$ can be written in the form
\[j_*[\PP(C\oplus\one)]=\sum_{i=0}^{n}e^G(L)^i(r^*\zeta_i)\]
for unique elements $\zeta_i\in\ChG(X)$.
Let $c$ be the rank of $C$ over $Z$, and $d=n-c$.  Then
{\allowdisplaybreaks
\bears
i_*s^G_k(C) &=& i_*l_*e^G(L')^{k+c}([\PP(C\oplus\one)]) \\
&=& r_*j_*(e^G(L')^{k+c}([\PP(C\oplus\one)])) \\
&=& r_*e^G(L)^{k+c}(j_*[\PP(C\oplus\one)])) \\
&=& r_*\sum_{i=0}^{n}e^G(L)^{k+c+i}r^*(\zeta_i) \\
&=& \sum_{i=0}^{n}s^G_{k+i-d}(E)\cdot\zeta_i
\eears
}
Therefore:
\bears
\sum_{k}i_*s^G_k(C)t^k &=& \sum_{k}\sum_{i}(s^G_{k+i-d}(E)t^{k+i})(\zeta_{i}t^{-i}) \\
&=& (\sum_{j}s^G_{j-d}(E)t^j)(\sum_{i}\zeta_{i}t^{-i})
\eears
and so
\[\sum_{i=0}^{n}\zeta_{i}t^{n-i}=(\sum_{j}c^G_j(E)t^j)(\sum_{k}i_*s^G_k(C)t^k)\]
i.e. 
\[\zeta_i=\sum_{j+k=n-i}\!\!c^G_j(E)i_*s^G_k(C).\]

Let $\phi:E\to\PP(E\oplus\one)$ be the open embedding.  Since the pull-back $\phi^*L$ is a trivial line bundle over $E$, we have, for any $\delta\in\Ch{X}$,
\[
(\phi^*\circ e^G(L)^i)(r^*\delta)=e^G(\phi^*L)^i(\phi^*r^*\delta)=e^G(\phi^*L)^i(q^*\delta)=
\begin{cases}
q^*\delta & \textrm{if}\,\,\,i=0 \\
0 & \textrm{if}\,\,\,i>0
\end{cases}
\]
Therefore,
\bears
[C] &=& \phi^*([j_*[\PP(C\oplus\one)]) \\
&=& \phi^*(\sum_{i=0}^{n}e^G(L)^{i}r^*(\!\!\sum_{j+k=n-i}\!\!c^G_j(E)i_*s^G_k(C))) \\
&=& q^*(\sum_{j+k=n}c^G_j(E)i_*s^G_k(C))\in\ChG(E)
\eears
and the result follows by the definition of $C^G_E$.
\end{proof}

\begin{ex}\label{Xsmooth}  Suppose $X$ is smooth of dimension $d$ and that $C=E=N$, the normal bundle of $X\to R^p(X)$.  Then by Proposition \ref{subcone}, $C^G_E=(c^G(N)s^G(N))_{d}=[X]\in\ChG(X)$ and so $C_E=\wt{C}_E=[X]\in\Ch{X}$.
\end{ex}

\section{Properties of the Steenrod operations}\label{steenrodprop}

\subsection{Basic properties}

\begin{thm}\label{id0} $S^X_k=0$ if $k<0$ and $S^X_0$ is the identity.
\end{thm}
\begin{proof}
We first consider $S_k^X([X])$ where $X$ is a variety of dimension $d$.  By dimension count, $S^X_k([X])\in\Chk{d-k}{X}$ is trivial if $k<0$.  To compute $S_0^X([X])$, we may replace $X$ by $X_{red}$, pass to an algebraic closure, and replace $X$ by a smooth open subscheme.  Let $N$ be the normal bundle of $X$ in $R^p(X)=X^p$.  Then by Corollary \ref{corcalc}, Proposition \ref{subcone} and Example \ref{Xsmooth}, 
\[S^X_0([X])=b_0(-T_X)(\wt{N}_N)=[X]\]

In general, let $i:Z\to X$ be a closed subvariety.  Then by Proposition \ref{steenrodclosedembedding} and the first part of the proof, $S^X_k([Z])=i_*S^Z_k([Z])$ is trivial for $k<0$ and is equal to $[Z]\in\Ch{X}$ if $k=0$.
\end{proof}

\begin{prop}\label{smoothformula} (Formula for a smooth cycle) Let $Z$ be 
a smooth closed subvariety of a scheme $X$.  Then
\[S^X([Z])=i_*b(-T_Z)([Z])\]
where $i:Z\to X$ is the closed embedding.
\end{prop}
\begin{proof}First, assume $X$ is smooth; let $d$ be the dimension of $X$.  Let $N$ be the normal bundle of $X\to R^p(X)$.  Then by Corollary \ref{corcalc}, Proposition \ref{subcone} and Example \ref{Xsmooth}:
\[S^X([X])=b(-T_X)(\wt{N}_N)=b(-T_X)([X])\]

In general, if $Z$ is a smooth closed subvariety of a (not necessarily smooth) scheme $X$, then by Proposition \ref{steenrodclosedembedding} and the first part of the proof, we have
\[S^X([Z])=i_*S^Z([Z])=i_*b(-T_Z)([Z])\] 
\end{proof}

\begin{thm}\label{externalproducts} Let $X$ and $Y$ be two schemes.  Then
\(S^{X\times Y}(\gamma\times\delta )=S^X(\gamma )\times S^Y(\delta )\) for
any $\gamma\in\Ch{X}$ and $\delta\in\Ch{Y}$.
\end{thm}
\begin{proof}
We may assume that $\gamma=[Z]$ and $\delta=[B]$ where $Z\subset X$ is a closed subvariety of dimension $n$ and $B\subset Y$ is a closed subvariety of dimension $m$.

Let $f:X\to W$ and $g:B\to M$ be closed embeddings of $X, Y$ into smooth schemes $W$ and $M$ of dimensions $r$ and $s$, respectively.  Let $C$ and $D$ be the normal cones of $X\to R^p(X)$ and $Y\to R^p(Y)$, respectively.  Let $E$ be the restriction to $X$ of the normal bundle of $W\to R^p(W)$, and similarly, $V$ the restriction to $Y$ of the normal bundle of $M\to R^p(M)$.  Note that $f\times g: X\times Y\to W\times M$ is a closed embedding of $X\times Y$ into a smooth scheme of dimension $r+s$, that $E\times V$ is the restriction to $X\times Y$ of the normal bundle of $W\times M\to R^p(W\times M)$, and that the normal cone of $Z\times B\to R^p(Z\times B)$ is $C\times D$.

By Corollary \ref{corcalc}, it remains to show that
\begin{multline*}
(-1)^{(r+s)+(n+m)}b(-T_{W\times M}|_{X\times Y})(\wt{(C\times D)}_{(E\times V)}) \\
=(-1)^{r+n}b(-T_W|_X)(\wt{C}_E)\times (-1)^{s+m}b(-T_M|_Y)(\wt{D}_V)
\end{multline*}
Since $b$ is multiplicative, it is enough to show that
\[\wt{(C\times D)}_{(E\times V)}=\wt{C}_E\times\wt{D}_V\]
but this follows from the fact that $[C]\times [D]=[C\times D]\in\ChG(E\times V).$
\end{proof}

\subsection{Functorialities}

\begin{lemma}\label{PrX} 
Let $q:\PP^r\times X\to X$ be the projection.  Then $q_*\circ 
S^{\PP^r\times X}=S^X\circ q_*$.
\end{lemma}
\begin{proof}
The group $\CH{\PP^r\times X}$ is generated by the cycles $\gamma 
=[\PP^k\times Z]$ for all closed subvarieties $Z\subset X$ and $k\leq r$.
Hence by Proposition \ref{steenrodclosedembedding} we may assume $Z=X$ and 
$k=r$.  The statement is obvious if $r=0$, so we may assume $r>0$.
Since $q_*(\gamma )=0$ we need to prove that $q_*S^{\PP^r\times X}(\gamma )
=0$.

By Theorem \ref{externalproducts}, we have that
\[S^{\PP^r\times X}(\gamma )=S^{\PP^r}([\PP^r])\times S^X([X])\]
but by Proposition \ref{smoothformula}, $S^{\PP^r}([\PP^r])=
b(-T_{\PP^r})([\PP^r])=(1+h^{p-1})^{-r-1}$, where $h$ is the class of a 
hyperplane in $\PP^r$.  Therefore, 
\[q_*S^{\PP^r\times X}(\gamma )=\deg (1+h^{p-1})^{-r-1}\cdot S^X([X])\]
Now $\deg (1+h^{p-1})^{-r-1}=0$ unless $r=(p-1)k$ for some integer $k$, 
and in this case, 
\[\deg (1+h^{p-1})^{-r-1}=\binom{-r-1}{k}=\binom{-(p-1)k-1}{k}
=(-1)^k\binom{pk}{k}\]
which is divisible by $p$ if $k>0$.
\end{proof}

\begin{thm}\label{steenrodproper} Let $f:X\to Y$ be a proper morphism.  Then the following diagram commutes:
\[
\xymatrix{
\textnormal{Ch}(X) \ar[r]^{S^X} \ar[d]^{f_*} & \Ch{X} \ar[d]^{f_*} \\
\Ch{Y} \ar[r]^{S^Y} & \Ch{Y}
}
\]
\end{thm}
\begin{proof}
$f$ factors as the composition of a closed embedding $X\to\PP^r\times Y$
and the projection $\PP^r\times Y\to Y$, so this result follows from
Proposition \ref{steenrodclosedembedding} and Lemma \ref{PrX}.
\end{proof}

\begin{prop}\label{mu-rho-yx} Let $f:Y\to X$ be a regular closed embedding of  schemes of codimension $s$.  Let $V\to Y$ and $W\to X$ be vector bundles with $V^G=Y$ and $W^G=X$.  
Let $k:\bu{V}\to\bu{W}$ be the induced closed embedding.  Then
\[\rho^V\circ k^*(\gamma)=\mu (W|_Y/V)\circ f^*\circ\rho^W(\gamma)\]
for all $\gamma\in\Ch{\bu{W}}$.  In particular, if both $X$ and $Y$ are smooth, and if $W$ is the normal bundle of $X\to R^p(X)$ and $V$ is the normal bundle of $Y\to R^p(Y)$,
\[\rho^V\circ k^*(\gamma)=(-1)^s\omega (N)\circ f^*\circ\rho^W(\gamma)\]
where $N$ is the normal bundle of $Y$ in $X$.
\end{prop}
\begin{proof}Consider the diagram
\[
\xymatrix{
\bu{V} \ar[d] \ar[r]^i & \bu{E} \ar[d] \ar[r]^j & \bu{W} \ar[d] \\
Y \ar@{}[r]|-{=} & Y \ar[r]^f & X
}
\]
where $E=W|_Y$.  Then, by Corollary \ref{transversediagram} and Lemma \ref{mu-rho},
\[\rho^V(k^*(\gamma))=\rho^V(i^*(j^*(\gamma)))
=(\mu(E/V)\circ\rho^E)(j^*(\gamma))=(\mu(E/V)\circ f^*\circ\rho^W)(\gamma)\]

Now consider the case where $W$ is the normal bundle of $X\to R^p(X)$ and $V$ is the normal bundle of $Y\to R^p(Y)$; we proceed as in the proof of Corollary \ref{omegacor}.  As before, we set $H$ to be the cokernel of the embedding $F\to F[G]$; then $W=T_X\otimes H$, $V=T_Y \otimes H$, and $W|_Y/V=N\otimes H$.  Note that $s$ is the rank of $N$ as a vector bundle over $Y$.  Suppose $N$ has a filtration with quotients $N_i$ ($i=1,2,\ldots s$) and $H$ has a filtration with quotients $H_j$ ($j=1,2,\ldots p-1$); note that the $H_j$ are trivial line bundles.  Hence $N\otimes H$ has a filtration with quotients $N_i\otimes H_j$, where $G$ acts with weight $j$, and so
{\allowdisplaybreaks
\bears
\mu(N\otimes H) &=& \prod_{j=1}^{p-1}\prod_{i=1}^{s}(j+c_1(N_i\otimes H_j)) \\
&=& \prod_{j=1}^{p-1}\prod_{i=1}^{s}(j+c_1(N_i)+c_1(H_j)) \\
&=& \prod_{i=1}^{s}\prod_{j=1}^{p-1}(j+c_1(N_i)) \\
&=& \prod_{i=1}^{s}(-1+c_1(N_i)^{p-1}) \\
&=& (-1)^{s}\prod_{i=1}^{s}(1-c_1(N_i)^{p-1}) \\
&=& (-1)^s\omega(N)
\eears
}
\end{proof}

\begin{prop}\label{bSteenrod}Let $f:Y\to X$ be a closed embedding of smooth schemes with normal bundle $N$.  Then
\begin{equation}\label{beqn}
(b(N)\circ f^*\circ S^X)(\gamma)=(S^Y\circ f^*)(\gamma)
\end{equation}
for all $\gamma\in\Ch{X}$.
\end{prop}
\begin{proof}
We may assume $\gamma=[Z]$ for some closed subvariety $Z\subset X$.  Let $n=\dim (Z)$, $s$ the codimension of $f$, and $r=n-s=\dim f^*Z$.  Let $k:\bu{Y}\to\bu{X}$ be the morphism induced by $f$.  Then, by Propositions \ref{Pfunctoriality} and \ref{mu-rho-yx}:
\bears
&& (\sum_{i+j=k}b_{i(p-1)}(N)\circ f^*\circ S_{j}^X)([Z]) \\
&=& (\sum_{i+j=k}(-1)^{n+j}b_{i(p-1)}(N)\circ f^*\circ\rho_{(p-1)(n+j)}^X\circ P_X)([Z]) \\
&=& (\sum_{i+j=k}(-1)^{n+j}(-1)^i\omega_{i(p-1)}(N)\circ f^*\circ\rho_{(p-1)(n+j)}^X\circ P_X)
([Z]) \\
&=& (\sum_{i+j=k}(-1)^{n+j}(-1)^i(-1)^s\rho^Y_{(p-1)(r+k)}\circ k^*\circ P_X)([Z]) \\
&=& (\sum_{i+j=k}(-1)^{n-s}(-1)^{i+j}\rho^Y_{(p-1)(r+k)}\circ P_Y\circ f^*)([Z]) \\
&=& (-1)^{r+k}(\rho^Y_{(p-1)(r+k)}\circ P_Y\circ f^*)([Z]) \\
&=& (S_k^Y\circ f^*)([Z])
\eears

\end{proof}

\section{Steenrod operations on smooth schemes}\label{smoothsteenrod}

\begin{defn} Let $X$ be a smooth scheme.  We define the \textit{Steenrod operations of
cohomological type} by the formula
\[S_X=b(T_X)\circ S^X\]
We write $S_X=\coprod S^k_X$, where for each $k$, $S^k_X$ is an operation
\[S^k_X:\textrm{Ch}^n(X)\to \textrm{Ch}^{n+k(p-1)}(X)\]
\end{defn}

\begin{prop}\label{wu}(Wu Formula) Let $X$ be a smooth scheme, $Z\subset X$ a smooth closed subscheme.  Then $S_X([Z])=i_*b(N)([Z])$, where $N$ is the normal bundle of the closed embedding $i:Z\to X$.
\end{prop}
\begin{proof}By Proposition \ref{smoothformula}, we have
\bears
S_X([Z]) &=& (b(T_X)\circ S^X)([Z]) \\
&=& (b(T_X)\circ i_*\circ b(-T_Z))([Z]) \\
&=& (i_* b(i^*T_X)\circ b(-T_Z))([Z]) \\
&=& i_*b(N)([Z])
\eears
\end{proof}

\begin{thm}\label{external}Let $X$ and $Y$ be two smooth schemes.  Then $S_{X\times Y}=S_X\times S_Y$.
\end{thm}
\begin{proof}
Since $T_{X\times Y}=T_X\times T_Y$, we have by Theorem \ref{externalproducts} that
\[S_{X\times Y}=b(T_{X\times Y})\circ S^{X\times Y}=b(T_X)S^X\times b(T_Y)S^Y=S_X\times S_Y\]
\end{proof}

\begin{thm}\label{steenrodpullbacks}
Let $f:Y\to X$ be a morphism of smooth schemes.  Then $f^*\circ S_X=S_Y\circ f^*$.
\end{thm}
\begin{proof}
Suppose first that $f$ is a closed embedding with normal bundle $N$.  Then by Proposition \ref{bSteenrod},
{\allowdisplaybreaks
\bears
f^*\circ S_X &=& f^*\circ b(T_X)\circ S^X \\
&=& b(f^*T_X)\circ f^*\circ S^X \\
&=& b(T_Y)\circ b(N)\circ f^*\circ S^X \\
&=& b(T_Y)\circ S^Y\circ f^* \\
&=& S_Y\circ f^*
\eears
}
Second, suppose $f$ is a projection $f:Y\times X\to X$.  Let $Z\subset X$ be a closed
subvariety.  We have that $f^*([Z])=[Y\times Z]=[Y]\times [Z]$, so, using Proposition \ref{smoothformula} and Theorem \ref{external}:
\bears
S_{Y\times X}(f^*([Z])) &=& S_{Y\times X}([Y\times Z]) = S_{Y\times X}([Y]\times [Z]) \\
&=& S_Y([Y])\times S_X([Z]) \\
&=& (b(T_Y)S^Y([Y]))\times S_X([Z]) \\
&=& [Y]\times S_X([Z]) \\
&=& f^*S_X([Z])
\eears
In general, we can write $f$ as the composition of the closed embedding $Y\to Y\times X$
and the projection $Y\times X\to X$, so the result follows from the two above cases.
\end{proof}

\begin{cor} (Cartan Formula) Let $X$ be a smooth scheme.  Then
\[S_X(\gamma\cdot\delta)=S_X(\gamma)\cdot S_X(\delta)\]
for any $\gamma,\delta\in\Ch{X}$.
\end{cor}
\begin{proof}Let $i:X\to X\times X$ be the diagonal embedding.  Then by Theorems \ref{external} and \ref{steenrodpullbacks},
\bears
S_X(\gamma\cdot\delta)=S_X(i^*(\gamma\times\delta))
&=& i^*(S_{X\times Y}(\gamma\times\delta) \\
&=& i^*(S_X(\gamma)\times S_X(\delta)) \\
&=& S_X(\gamma)\cdot S_X(\delta)
\eears
\end{proof}

\begin{thm}\label{pthpower}Let $X$ be a smooth scheme.  Then for any 
$\delta\in\textnormal{Ch}^k(X)$,
\[
S^r_X(\delta)=
\begin{cases}
\delta, & r=0 \\
\delta^p, & r=k \\
0, & r<0\,\,\textnormal{or}\,\, r>k
\end{cases}
\]
\end{thm}
\begin{proof}By definition and Theorem \ref{id0}, $S^r_X=0$ if $r<0$ and $S_X^0$ is the
identity function.  We assume that $X$ is a variety and that $\delta=[Z]$ for a closed subvariety $Z\subset X$.  Let $d=\dim X, n=\dim Z$ so that $k=d-n$.  In the notation of Corollary \ref{corcalc}, we have
\[S_X(\delta)=b(T_X)b(-T_X)(\wt{C}_E)=\wt{C}_E=\sum\tilde{\gamma}_i\]
Since $\tilde{\gamma}_i\in\Chk{n-i(p-1)}{X}, S_X^r(\delta)=\tilde{\gamma}_r.$  Following the notation of Lemma \ref{brolemma}, $\tilde{\gamma}_r=(-1)^{d+n+r}a_{(d-n-r)(p-1)}$, and since $a_j=0$ for $j<0$, $S_X^r(\delta)=0$ if $r>d-n=k$.

It remains to show that $S_X^k(\delta)=\delta^p$.  Since $S_X^k(\delta)=\tilde{\gamma}_k=a_0$, we need to show that $\delta^p=a_0$.  By a restriction-corestriction argument, we may assume that $F$ contains $p$-th roots of unity, so that $R^p(X)=X^p$.  Let $\Delta:X\to X^p$ be the diagonal embedding, $C$ the normal cone of $Z\to Z^p$, and $N$ the normal bundle of $X\to X^p$.  Let $h:C\to N$ be the closed embedding, and $s:N\to X$ the projection.  By \cite{merkurjev}, Proposition 52.6,
\[s^*(\delta^p)=s^*\Delta^*([Z^p])=\sigma_{\Delta}([Z^p])=h_*([C])\]
where $\sigma_{\Delta}$ is the deformation homomorphism.  Now consider the diagram
\begin{equation}\label{lasteq}
\xymatrix{
\ChG(N) \ar[r] & \Ch{N} \\
\ChG(X) \ar[u]^{s_G^*} \ar[r] & \Ch{X} \ar[u]^{s^*}
}
\end{equation}
where the horizontal maps ``forget the $G$-action'' on $N$ and $X$ (to be precise, these maps are the transfer from $G$ to the trivial subgroup; see \cite{brosnan}, 3.4 for a detailed explanation).  We have that (\ref{lasteq}) commutes, as the horizontal maps are induced by a push-forward.  Consider $C^G_N\in\ChG(X)$; if we follow it right and up, we get $s^*(a_0)$; if we follow it up and then right, we get $h_*([C])$.  Hence $s^*(a_0)=h_*([C])=s^*(\delta^p)$, and therefore $\delta^p=a_0$, as $s^*$ is an isomorphism.
\end{proof}


\nocite{*}

\begin{thebibliography}{1}

\bibitem[1]{brosnan}
Patrick Brosnan.
\newblock ``Steenrod operations in {C}how theory.''
\newblock {\em Trans. Amer. Math. Soc.}, {\bfseries 355}(5):1869--1903
  (electronic), 2003.

\bibitem[2]{edidin}
Dan Edidin and William Graham.
\newblock ``Equivariant intersection theory.''
\newblock {\em Invent. Math.}, {\bfseries 131}(3):595--634, 1998.

\bibitem[3]{eisenbud}
David Eisenbud and Joe Harris.
\newblock {\em The geometry of schemes}, volume 197 of {\em Graduate Texts in
  Mathematics}.
\newblock Springer-Verlag, New York, 2000.

\bibitem[4]{merkurjev}
Richard Elman, Nikita Karpenko, and Alexander Merkurjev.
\newblock {\em The Algebraic and Geometric Theory of Quadratic Forms}.
\newblock Preprint.

\bibitem[5]{fulton}
William Fulton.
\newblock {\em Intersection theory}, volume~2 of {\em Ergebnisse der Mathematik
  und ihrer Grenzgebiete. 3. Folge. A Series of Modern Surveys in Mathematics
  [Results in Mathematics and Related Areas. 3rd Series. A Series of Modern
  Surveys in Mathematics]}.
\newblock Springer-Verlag, Berlin, second edition, 1998.

\bibitem[6]{hartshorne}
Robin Hartshorne.
\newblock {\em Algebraic geometry}.
\newblock Springer-Verlag, New York, 1977.
\newblock Graduate Texts in Mathematics, No. 52.

\bibitem[7]{thebook}
Max-Albert Knus, Alexander Merkurjev, Markus Rost, and Jean-Pierre Tignol.
\newblock {\em The book of involutions}, volume~44 of {\em American
  Mathematical Society Colloquium Publications}.
\newblock American Mathematical Society, Providence, RI, 1998.
\newblock With a preface in French by J.\ Tits.

\bibitem[8]{milne}
James~S. Milne.
\newblock {\em \'{E}tale cohomology}, volume~33 of {\em Princeton Mathematical
  Series}.
\newblock Princeton University Press, Princeton, N.J., 1980.

\bibitem[9]{rost}
Markus Rost.
\newblock ``Notes on the degree formula.''\\
\newblock http://www.mathematik.uni-bielefeld.de/\~{}rost/data/bd.pdf, 2001.

\bibitem[10]{steenrod}
N.~E. Steenrod.
\newblock {\em Cohomology operations}.
\newblock Lectures by N. E. Steenrod written and revised by D. B. A. Epstein.
  Annals of Mathematics Studies, No. 50. Princeton University Press, Princeton,
  N.J., 1962.

\bibitem[11]{voevodsky}
Vladimir Voevodsky.
\newblock ``The Milnor Conjecture.''
\newblock Preprint, December 20, 1996, K-theory Preprint Archives,
  http://www.math.uiuc.edu/K-theory/0170/.

\end{thebibliography}
\end{document}